\documentclass[12pt,final,notitlepage,onecolumn]{article}
\usepackage{amsfonts}
\usepackage{amssymb}
\usepackage{graphicx}
\usepackage{amsmath}

\input{tcilatex}

\begin{document}

\title{Generalizations of Popoviciu's inequality}
\author{Darij Grinberg}
\date{20 March 2008}
\maketitle

\begin{abstract}
We establish a general criterion for inequalities of the kind%
\begin{align*}
& \text{convex combination of }f\left( x_{1}\right) ,\text{ }f\left(
x_{2}\right) ,\text{ }...,\text{ }f\left( x_{n}\right) \\
& \ \ \ \ \ \ \ \ \ \ \text{ and }f\left( \text{some weighted mean of }x_{1},%
\text{ }x_{2},\text{ }...,\text{ }x_{n}\right) \\
& \geq \text{convex combination of }f\left( \text{some other weighted means
of }x_{1},\text{ }x_{2},\text{ }...,\text{ }x_{n}\right) ,
\end{align*}%
where $f$ is a convex function on an interval $I\subseteq \mathbb{R}$
containing the reals $x_{1},$ $x_{2},$ $...,$ $x_{n},$ to hold. Here, the
left hand side contains only one weighted mean, while the right hand side
may contain as many as possible, as long as there are finitely many. The
weighted mean on the left hand side must have positive weights, while those
on the right hand side must have nonnegative weights.

This criterion entails Vasile C\^{\i}rtoaje's generalization of the
Popoviciu inequality (in its standard and in its weighted forms) as well as
a cyclic inequality that sharpens another result by Vasile C\^{\i}rtoaje.
The latter cyclic inequality (in its non-weighted form) states that%
\begin{equation*}
2\sum_{i=1}^{n}f\left( x_{i}\right) +n\left( n-2\right) f\left( x\right)
\geq n\sum_{s=1}^{n}f\left( x+\dfrac{x_{s}-x_{s+r}}{n}\right) ,
\end{equation*}%
where indices are cyclic modulo $n,$ and $x=\dfrac{x_{1}+x_{2}+...+x_{n}}{n}%
. $
\end{abstract}

{\small This is the standard version of this note. A "formal" version with
more detailed proofs can be found at\newline
\texttt{http://www.stud.uni-muenchen.de/\symbol{126}%
darij.grinberg/PopoviciuFormal.pdf}\newline
However, due to these details, it is longer and much more troublesome to
read, so it should be used merely as a resort in case you do not understand
the proofs in this standard version.}

\bigskip


\begin{center}
\textbf{1. Introduction}
\end{center}

The last few years saw some activity related to the Popoviciu inequality on
convex functions. Some generalizations were conjectured and subsequently
proven using majorization theory and (mostly) a lot of computations. In this
note I am presenting an apparently new approach that proves these
generalizations as well as some additional facts with a lesser amount of
computation and avoiding majorization theory (more exactly, avoiding the
standard, asymmetric definition of majorization; we will prove a "symmetric"
version of the Karamata inequality on the way, which will not even use the
word "majorize").

The very starting point of the whole theory is the following famous fact:

\begin{quote}
\textbf{Theorem 1a, the Jensen inequality.} Let $f$ be a convex function
from an interval $I\subseteq\mathbb{R}$ to $\mathbb{R}.$ Let $x_{1},$ $%
x_{2}, $ $...,$ $x_{n}$ be finitely many points from $I.$ Then,%
\begin{equation*}
\dfrac{f\left( x_{1}\right) +f\left( x_{2}\right) +...+f\left( x_{n}\right) 
}{n}\geq f\left( \dfrac{x_{1}+x_{2}+...+x_{n}}{n}\right) .
\end{equation*}

In words, the arithmetic mean of the values of $f$ at the points $x_{1},$ $%
x_{2},$ $...,$ $x_{n}$ is greater or equal to the value of $f$ at the
arithmetic mean of these points.
\end{quote}

We can obtain a "weighted version" of this inequality by replacing
arithmetic means by weighted means with some nonnegative weights $w_{1},$ $%
w_{2},$ $...,$ $w_{n}$:

\begin{quote}
\textbf{Theorem 1b, the weighted Jensen inequality.} Let $f$ be a convex
function from an interval $I\subseteq\mathbb{R}$ to $\mathbb{R}.$ Let $%
x_{1}, $ $x_{2},$ $...,$ $x_{n}$ be finitely many points from $I.$ Let $%
w_{1}, $ $w_{2},$ $...,$ $w_{n}$ be $n$ nonnegative reals which are not all
equal to $0.$ Then,%
\begin{equation*}
\dfrac{w_{1}f\left( x_{1}\right) +w_{2}f\left( x_{2}\right)
+...+w_{n}f\left( x_{n}\right) }{w_{1}+w_{2}+...+w_{n}}\geq f\left( \dfrac{%
w_{1}x_{1}+w_{2}x_{2}+...+w_{n}x_{n}}{w_{1}+w_{2}+...+w_{n}}\right) .
\end{equation*}
\end{quote}

Obviously, Theorem 1a follows from Theorem 1b applied to $%
w_{1}=w_{2}=...=w_{n}=1,$ so that Theorem 1b is more general than Theorem 1a.

We won't stop at discussing equality cases here, since they can depend in
various ways on the input (i. e., on the function $f,$ the reals $w_{1},$ $%
w_{2},$ $...,$ $w_{n}$ and the points $x_{1},$ $x_{2},$ $...,$ $x_{n}$) -
but each time we use a result like Theorem 1b, with enough patience we can
extract the equality case from the proof of this result and the properties
of the input.

The Jensen inequality, in both of its versions above, is applied often
enough to be called one of the main methods of proving inequalities. Now, in
1965, a similarly styled inequality was found by the Romanian Tiberiu
Popoviciu:

\begin{quote}
\textbf{Theorem 2a, the Popoviciu inequality.} Let $f$ be a convex function
from an interval $I\subseteq\mathbb{R}$ to $\mathbb{R},$ and let $x_{1},$ $%
x_{2},$ $x_{3}$ be three points from $I.$ Then,%
\begin{equation*}
f\left( x_{1}\right) +f\left( x_{2}\right) +f\left( x_{3}\right) +3f\left( 
\dfrac{x_{1}+x_{2}+x_{3}}{3}\right) \geq2f\left( \dfrac {x_{2}+x_{3}}{2}%
\right) +2f\left( \dfrac{x_{3}+x_{1}}{2}\right) +2f\left( \dfrac{x_{1}+x_{2}%
}{2}\right) .
\end{equation*}
\end{quote}

Again, a weighted version can be constructed:

\begin{quote}
\textbf{Theorem 2b, the weighted Popoviciu inequality.} Let $f$ be a convex
function from an interval $I\subseteq\mathbb{R}$ to $\mathbb{R},$ let $%
x_{1}, $ $x_{2},$ $x_{3}$ be three points from $I,$ and let $w_{1},$ $w_{2},$
$w_{3} $ be three nonnegative reals such that $w_{2}+w_{3}\neq0,$ $%
w_{3}+w_{1}\neq0$ and $w_{1}+w_{2}\neq0.$ Then,%
\begin{align*}
& w_{1}f\left( x_{1}\right) +w_{2}f\left( x_{2}\right) +w_{3}f\left(
x_{3}\right) +\left( w_{1}+w_{2}+w_{3}\right) f\left( \dfrac{%
w_{1}x_{1}+w_{2}x_{2}+w_{3}x_{3}}{w_{1}+w_{2}+w_{3}}\right) \\
& \geq\left( w_{2}+w_{3}\right) f\left( \dfrac{w_{2}x_{2}+w_{3}x_{3}}{%
w_{2}+w_{3}}\right) +\left( w_{3}+w_{1}\right) f\left( \dfrac{%
w_{3}x_{3}+w_{1}x_{1}}{w_{3}+w_{1}}\right) +\left( w_{1}+w_{2}\right)
f\left( \dfrac{w_{1}x_{1}+w_{2}x_{2}}{w_{1}+w_{2}}\right) .
\end{align*}
\end{quote}

The really interesting part of the story began when Vasile C\^{\i}rtoaje -
alias "Vasc" on the MathLinks forum - proposed the following two
generalizations of Theorem 2a ([1] and [2] for Theorem 3a, and [1] and [3]
for Theorem 4a):

\begin{quote}
\textbf{Theorem 3a (Vasile C\^{\i}rtoaje).} Let $f$ be a convex function
from an interval $I\subseteq\mathbb{R}$ to $\mathbb{R}.$ Let $x_{1},$ $%
x_{2}, $ $...,$ $x_{n}$ be finitely many points from $I.$ Then,%
\begin{equation*}
\sum_{i=1}^{n}f\left( x_{i}\right) +n\left( n-2\right) f\left( \dfrac{%
x_{1}+x_{2}+...+x_{n}}{n}\right) \geq\sum_{j=1}^{n}\left( n-1\right) f\left( 
\dfrac{\sum\limits_{1\leq i\leq n;\ i\neq j}x_{i}}{n-1}\right) .
\end{equation*}

\textbf{Theorem 4a (Vasile C\^{\i}rtoaje).} Let $f$ be a convex function
from an interval $I\subseteq\mathbb{R}$ to $\mathbb{R}.$ Let $x_{1},$ $%
x_{2}, $ $...,$ $x_{n}$ be finitely many points from $I.$ Then,%
\begin{equation*}
\left( n-2\right) \sum_{i=1}^{n}f\left( x_{i}\right) +nf\left( \dfrac{%
x_{1}+x_{2}+...+x_{n}}{n}\right) \geq\sum_{1\leq i<j\leq n}2f\left( \dfrac{%
x_{i}+x_{j}}{2}\right) .
\end{equation*}
\end{quote}

In [1], both of these facts were nicely proven by C\^{\i}rtoaje. I gave a
different and rather long proof of Theorem 3a in [2]. All of these proofs
use the Karamata inequality. Theorem 2a follows from each of the Theorems 3a
and 4a upon setting $n=3.$

It is pretty straightforward to obtain generalizations of Theorems 3a and 4a
by putting in weights as in Theorems 1b and 2b. A more substantial
generalization was given by Yufei Zhao - alias "Billzhao" on MathLinks - in
[3]:

\begin{quote}
\textbf{Theorem 5a (Yufei Zhao).} Let $f$ be a convex function from an
interval $I\subseteq\mathbb{R}$ to $\mathbb{R}.$ Let $x_{1},$ $x_{2},$ $..., 
$ $x_{n}$ be finitely many points from $I,$ and let $m$ be an integer. Then,%
\begin{align*}
& \binom{n-2}{m-1}\sum_{i=1}^{n}f\left( x_{i}\right) +\binom{n-2}{m-2}%
nf\left( \frac{x_{1}+x_{2}+...+x_{n}}{n}\right) \\
& \geq\sum_{1\leq i_{1}<i_{2}<...<i_{m}\leq n}mf\left( \frac{%
x_{i_{1}}+x_{i_{2}}+...+x_{i_{m}}}{m}\right) .
\end{align*}
\end{quote}

Note that if $m\leq0$ or $m>n,$ the sum $\sum\limits_{1\leq
i_{1}<i_{2}<...<i_{m}\leq n}mf\left( \dfrac{x_{i_{1}}+x_{i_{2}}+...+x_{i_{m}}%
}{m}\right) $ is empty, so that its value is $0.$

Note that Theorems 3a and 4a both are particular cases of Theorem 5a (in
fact, set $m=n-1$ to get Theorem 3a and $m=2$ to get Theorem 4a).

An rather complicated proof of Theorem 5a was given by myself in [3]. After
some time, the MathLinks user "Zhaobin" proposed a weighted version of this
result:

\begin{quote}
\textbf{Theorem 5b (Zhaobin).} Let $f$ be a convex function from an interval 
$I\subseteq\mathbb{R}$ to $\mathbb{R}.$ Let $x_{1},$ $x_{2},$ $...,$ $x_{n}$
be finitely many points from $I,$ let $w_{1},$ $w_{2},$ $...,$ $w_{n}$ be
nonnegative reals, and let $m$ be an integer. Assume that $%
w_{1}+w_{2}+...+w_{n}\neq0,$ and that $w_{i_{1}}+w_{i_{2}}+...+w_{i_{m}}%
\neq0 $ for any $m$ integers $i_{1},$ $i_{2},$ $...,$ $i_{m}$ satisfying $%
1\leq i_{1}<i_{2}<...<i_{m}\leq n.$

Then,%
\begin{align*}
& \binom{n-2}{m-1}\sum_{i=1}^{n}w_{i}f\left( x_{i}\right) +\binom{n-2}{m-2}%
\left( w_{1}+w_{2}+...+w_{n}\right) f\left( \frac{%
w_{1}x_{1}+w_{2}x_{2}+...+w_{n}x_{n}}{w_{1}+w_{2}+...+w_{n}}\right) \\
& \geq\sum_{1\leq i_{1}<i_{2}<...<i_{m}\leq n}\left(
w_{i_{1}}+w_{i_{2}}+...+w_{i_{m}}\right) f\left( \frac{%
w_{i_{1}}x_{i_{1}}+w_{i_{2}}x_{i_{2}}+...+w_{i_{m}}x_{i_{m}}}{%
w_{i_{1}}+w_{i_{2}}+...+w_{i_{m}}}\right) .
\end{align*}
\end{quote}

If we set $w_{1}=w_{2}=...=w_{n}=1$ in Theorem 5b, we obtain Theorem 5a. On
the other hand, putting $n=3$ and $m=2$ in Theorem 5b, we get Theorem 2b.

In this note, I am going to prove Theorem 5b (and therefore also its
particular cases - Theorems 2a, 2b, 3a, 4a and 5a). The proof is going to
use no preknowledge - in particular, classical majorization theory will be
avoided. Then, we are going to discuss an assertion similar to Theorem 5b
with its applications.

\begin{center}
\textbf{2. Absolute values interpolate convex functions}
\end{center}

We start preparing for our proof by showing a property of convex functions
which is definitely not new - it was mentioned by a MathLinks user called
"Fleeting\_Guest" in [4], post \#18 as a known fact:

\begin{quote}
\textbf{Theorem 6.} Let $f$ be a convex function from an interval $I\subseteq%
\mathbb{R}$ to $\mathbb{R}.$ Let $x_{1},$ $x_{2},$ $...,$ $x_{n}$ be
finitely many points from $I.$ Then, there exist two real constants $u$ and $%
v$ and $n$ nonnegative constants $a_{1},$ $a_{2},$ $...,$ $a_{n}$ such that%
\begin{equation*}
f\left( t\right) =vt+u+\sum\limits_{i=1}^{n}a_{i}\left\vert
t-x_{i}\right\vert \text{ holds for every }t\in\left\{
x_{1},x_{2},...,x_{n}\right\} .
\end{equation*}
\end{quote}

In brief, this result states that every convex function $f\left( x\right) $
on $n$ reals $x_{1},$ $x_{2},$ $...,$ $x_{n}$ can be interpolated by a
linear combination with nonnegative coefficients of a linear function and
the $n$ functions $\left\vert x-x_{i}\right\vert .$

The \textit{proof of Theorem 6}, albeit technical, will be given here for
the sake of completeness: First, we need an almost trivial fact which I use
to call the $\max\left\{ 0,x\right\} $ \textit{formula:} For any real number 
$x,$ we have $\max\left\{ 0,x\right\} =\dfrac{1}{2}\left( x+\left\vert
x\right\vert \right) .$

Furthermore, we denote $f\left[ y,z\right] =\dfrac{f\left( y\right) -f\left(
z\right) }{y-z}$ for any two points $y$ and $z$ from $I$ satisfying $y\neq
z. $ Then, we have $\left( y-z\right) \cdot f\left[ y,z\right] =f\left(
y\right) -f\left( z\right) $ for any two points $y$ and $z$ from $I$
satisfying $y\neq z.$

We can assume that all points $x_{1},$ $x_{2},$ $...,$ $x_{n}$ are
pairwisely distinct (if not, we can remove all superfluous $x_{i}$ and apply
Theorem 6 to the remaining points). Therefore, we can WLOG assume that $%
x_{1}<x_{2}<...<x_{n}.$ Then, for every $j\in\left\{ 1,2,...,n\right\} ,$ we
have%
\begin{align*}
f\left( x_{j}\right) & =f\left( x_{1}\right) +\sum_{k=1}^{j-1}\left( f\left(
x_{k+1}\right) -f\left( x_{k}\right) \right) =f\left( x_{1}\right)
+\sum_{k=1}^{j-1}\left( x_{k+1}-x_{k}\right) \cdot f\left[ x_{k+1},x_{k}%
\right] \\
& =f\left( x_{1}\right) +\sum_{k=1}^{j-1}\left( x_{k+1}-x_{k}\right)
\cdot\left( f\left[ x_{2},x_{1}\right] +\sum_{i=2}^{k}\left( f\left[
x_{i+1},x_{i}\right] -f\left[ x_{i},x_{i-1}\right] \right) \right) \\
& =f\left( x_{1}\right) +\sum_{k=1}^{j-1}\left( x_{k+1}-x_{k}\right) \cdot f 
\left[ x_{2},x_{1}\right] +\sum_{k=1}^{j-1}\left( x_{k+1}-x_{k}\right)
\cdot\sum_{i=2}^{k}\left( f\left[ x_{i+1},x_{i}\right] -f\left[ x_{i},x_{i-1}%
\right] \right) \\
& =f\left( x_{1}\right) +f\left[ x_{2},x_{1}\right] \cdot\sum_{k=1}^{j-1}%
\left( x_{k+1}-x_{k}\right) +\sum_{k=1}^{j-1}\sum_{i=2}^{k}\left( f\left[
x_{i+1},x_{i}\right] -f\left[ x_{i},x_{i-1}\right] \right) \cdot\left(
x_{k+1}-x_{k}\right) \\
& =f\left( x_{1}\right) +f\left[ x_{2},x_{1}\right] \cdot\sum_{k=1}^{j-1}%
\left( x_{k+1}-x_{k}\right) +\sum_{i=2}^{j-1}\sum_{k=i}^{j-1}\left( f\left[
x_{i+1},x_{i}\right] -f\left[ x_{i},x_{i-1}\right] \right) \cdot\left(
x_{k+1}-x_{k}\right) \\
& =f\left( x_{1}\right) +f\left[ x_{2},x_{1}\right] \cdot\sum_{k=1}^{j-1}%
\left( x_{k+1}-x_{k}\right) +\sum_{i=2}^{j-1}\left( f\left[ x_{i+1},x_{i}%
\right] -f\left[ x_{i},x_{i-1}\right] \right) \cdot\sum _{k=i}^{j-1}\left(
x_{k+1}-x_{k}\right) \\
& =f\left( x_{1}\right) +f\left[ x_{2},x_{1}\right] \cdot\left(
x_{j}-x_{1}\right) +\sum_{i=2}^{j-1}\left( f\left[ x_{i+1},x_{i}\right] -f%
\left[ x_{i},x_{i-1}\right] \right) \cdot\left( x_{j}-x_{i}\right) .
\end{align*}
Now we set%
\begin{align*}
& \alpha_{1}=\alpha_{n}=0; \\
& \alpha_{i}=f\left[ x_{i+1},x_{i}\right] -f\left[ x_{i},x_{i-1}\right] \ \
\ \ \ \ \ \ \ \ \text{for all }i\in\left\{ 2,3,...,n-1\right\} .
\end{align*}
Using these notations, the above computation becomes%
\begin{align*}
f\left( x_{j}\right) & =f\left( x_{1}\right) +f\left[ x_{2},x_{1}\right]
\cdot\left( x_{j}-x_{1}\right) +\sum_{i=2}^{j-1}\alpha _{i}\cdot\left(
x_{j}-x_{i}\right) \\
& =f\left( x_{1}\right) +f\left[ x_{2},x_{1}\right] \cdot\left(
x_{j}-x_{1}\right) +\underbrace{0}_{=\alpha_{1}}\cdot\max\left\{
0,x_{j}-x_{1}\right\} \\
& +\sum_{i=2}^{j-1}\alpha_{i}\cdot\underbrace{\left( x_{j}-x_{i}\right) }%
_{=\max\left\{ 0,x_{j}-x_{i}\right\} ,\text{ since }x_{j}-x_{i}\geq0,\text{
as }x_{i}\leq x_{j}}+\sum_{i=j}^{n}\alpha_{i}\cdot\underbrace{0}%
_{=\max\left\{ 0,x_{j}-x_{i}\right\} ,\text{ since }x_{j}-x_{i}\leq0,\text{
as }x_{j}\leq x_{i}} \\
& =f\left( x_{1}\right) +f\left[ x_{2},x_{1}\right] \cdot\left(
x_{j}-x_{1}\right) +\alpha_{1}\cdot\max\left\{ 0,x_{j}-x_{1}\right\} \\
& +\sum_{i=2}^{j-1}\alpha_{i}\cdot\max\left\{ 0,x_{j}-x_{i}\right\}
+\sum_{i=j}^{n}\alpha_{i}\cdot\max\left\{ 0,x_{j}-x_{i}\right\} \\
& =f\left( x_{1}\right) +f\left[ x_{2},x_{1}\right] \cdot\left(
x_{j}-x_{1}\right) +\sum_{i=1}^{n}\alpha_{i}\cdot\max\left\{
0,x_{j}-x_{i}\right\} \\
& =f\left( x_{1}\right) +f\left[ x_{2},x_{1}\right] \cdot\left(
x_{j}-x_{1}\right) +\sum_{i=1}^{n}\alpha_{i}\cdot\dfrac{1}{2}\left( \left(
x_{j}-x_{i}\right) +\left\vert x_{j}-x_{i}\right\vert \right) \\
& \ \ \ \ \ \ \ \ \ \ \left( \text{since }\max\left\{ 0,x_{j}-x_{i}\right\} =%
\frac{1}{2}\left( \left( x_{j}-x_{i}\right) +\left\vert
x_{j}-x_{i}\right\vert \right) \text{ by the }\max\left\{ 0,x\right\} \text{
formula}\right) \\
& =f\left( x_{1}\right) +f\left[ x_{2},x_{1}\right] \cdot\left(
x_{j}-x_{1}\right) +\sum_{i=1}^{n}\alpha_{i}\cdot\dfrac{1}{2}\left(
x_{j}-x_{i}\right) +\sum_{i=1}^{n}\alpha_{i}\cdot\dfrac{1}{2}\left\vert
x_{j}-x_{i}\right\vert \\
& =f\left( x_{1}\right) +\left( f\left[ x_{2},x_{1}\right] x_{j}-f\left[
x_{2},x_{1}\right] x_{1}\right) +\left( \dfrac{1}{2}\sum_{i=1}^{n}%
\alpha_{i}x_{j}-\dfrac{1}{2}\sum_{i=1}^{n}\alpha_{i}x_{i}\right)
+\sum_{i=1}^{n}\dfrac{1}{2}\alpha_{i}\left\vert x_{j}-x_{i}\right\vert \\
& =\left( f\left[ x_{2},x_{1}\right] +\dfrac{1}{2}\sum_{i=1}^{n}\alpha
_{i}\right) x_{j}+\left( f\left( x_{1}\right) -f\left[ x_{2},x_{1}\right]
x_{1}-\dfrac{1}{2}\sum_{i=1}^{n}\alpha_{i}x_{i}\right) +\sum_{i=1}^{n}\dfrac{%
1}{2}\alpha_{i}\left\vert x_{j}-x_{i}\right\vert .
\end{align*}
Thus, if we denote%
\begin{align*}
v & =f\left[ x_{2},x_{1}\right] +\dfrac{1}{2}\sum\limits_{i=1}^{n}%
\alpha_{i};\ \ \ \ \ \ \ \ \ \ u=f\left( x_{1}\right) -f\left[ x_{2},x_{1}%
\right] x_{1}-\dfrac{1}{2}\sum_{i=1}^{n}\alpha_{i}x_{i}; \\
a_{i} & =\dfrac{1}{2}\alpha_{i}\ \ \ \ \ \ \ \ \ \ \text{for all }%
i\in\left\{ 1,2,...,n\right\} ,
\end{align*}
then we have%
\begin{equation*}
f\left( x_{j}\right) =vx_{j}+u+\sum\limits_{i=1}^{n}a_{i}\left\vert
x_{j}-x_{i}\right\vert .
\end{equation*}
Since we have shown this for every $j\in\left\{ 1,2,...,n\right\} ,$ we can
restate this as follows: We have%
\begin{equation*}
f\left( t\right) =vt+u+\sum\limits_{i=1}^{n}a_{i}\left\vert
t-x_{i}\right\vert \text{ for every }t\in\left\{
x_{1},x_{2},...,x_{n}\right\} .
\end{equation*}
Hence, in order for the proof of Theorem 6 to be complete, it is enough to
show that the $n$ reals $a_{1},$ $a_{2},$ $...,$ $a_{n}$ are nonnegative.
Since $a_{i}=\dfrac{1}{2}\alpha_{i}$ for every $i\in\left\{
1,2,...,n\right\} ,$ this will follow once it is proven that the $n$ reals $%
\alpha_{1},$ $\alpha_{2},$ $...,$ $\alpha_{n}$ are nonnegative. Thus, we
have to show that $\alpha_{i}$ is nonnegative for every $i\in\left\{
1,2,...,n\right\} .$ This is trivial for $i=1$ and for $i=n$ (since $%
\alpha_{1}=0$ and $\alpha_{n}=0$), so it remains to prove that $\alpha_{i}$
is nonnegative for every $i\in\left\{ 2,3,...,n-1\right\} .$ Now, since $%
\alpha_{i}=f\left[ x_{i+1},x_{i}\right] -f\left[ x_{i},x_{i-1}\right] $ for
every $i\in\left\{ 2,3,...,n-1\right\} ,$ we thus have to show that $f\left[
x_{i+1},x_{i}\right] -f\left[ x_{i},x_{i-1}\right] $ is nonnegative for
every $i\in\left\{ 2,3,...,n-1\right\} .$ In other words, we have to prove
that $f\left[ x_{i+1},x_{i}\right] \geq f\left[ x_{i},x_{i-1}\right] $ for
every $i\in\left\{ 2,3,...,n-1\right\} .$ But since $x_{i-1}<x_{i}<x_{i+1},$
this follows from the next lemma:

\begin{quote}
\textbf{Lemma 7.} Let $f$ be a convex function from an interval $I\subseteq 
\mathbb{R}$ to $\mathbb{R}.$ Let $x,$ $y,$ $z$ be three points from $I$
satisfying $x<y<z.$ Then, $f\left[ z,y\right] \geq f\left[ y,x\right] .$
\end{quote}

\textit{Proof of Lemma 7.} Since the function $f$ is convex on $I,$ and
since $z$ and $x$ are points from $I,$ the definition of convexity yields 
\begin{equation*}
\dfrac{\dfrac{1}{z-y}f\left( z\right) +\dfrac{1}{y-x}f\left( x\right) }{%
\dfrac{1}{z-y}+\dfrac{1}{y-x}}\geq f\left( \dfrac{\dfrac{1}{z-y}z+\dfrac {1}{%
y-x}x}{\dfrac{1}{z-y}+\dfrac{1}{y-x}}\right)
\end{equation*}
(here we have used that $\dfrac{1}{z-y}>0$ and $\dfrac{1}{y-x}>0,$ what is
clear from $x<y<z$). Since $\dfrac{\dfrac{1}{z-y}z+\dfrac{1}{y-x}x}{\dfrac {1%
}{z-y}+\dfrac{1}{y-x}}=y,$ this simplifies to%
\begin{align*}
\dfrac{\dfrac{1}{z-y}f\left( z\right) +\dfrac{1}{y-x}f\left( x\right) }{%
\dfrac{1}{z-y}+\dfrac{1}{y-x}} & \geq f\left( y\right) ,\ \ \ \ \ \ \ \ \ \ 
\text{so that} \\
\dfrac{1}{z-y}f\left( z\right) +\dfrac{1}{y-x}f\left( x\right) & \geq\left( 
\dfrac{1}{z-y}+\dfrac{1}{y-x}\right) f\left( y\right) ,\ \ \ \ \ \ \ \ \ \ 
\text{so that} \\
\dfrac{1}{z-y}f\left( z\right) -\dfrac{1}{z-y}f\left( y\right) & \geq\dfrac{1%
}{y-x}f\left( y\right) -\dfrac{1}{y-x}f\left( x\right) ,\ \ \ \ \ \ \ \ \ \ 
\text{so that} \\
\dfrac{f\left( z\right) -f\left( y\right) }{z-y} & \geq\dfrac{f\left(
y\right) -f\left( x\right) }{y-x}.
\end{align*}
This becomes $f\left[ z,y\right] \geq f\left[ y,x\right] ,$ and thus Lemma 7
is proven. Thus, the proof of Theorem 6 is completed.

\begin{center}
\textbf{3. The Karamata inequality in symmetric form}
\end{center}

Now as Theorem 6 is proven, it becomes easy to prove the Karamata inequality
in the following form:

\begin{quote}
\textbf{Theorem 8a, the Karamata inequality in symmetric form.} Let $f$ be a
convex function from an interval $I\subseteq\mathbb{R}$ to $\mathbb{R},$ and
let $n$ be a positive integer. Let $x_{1},$ $x_{2},$ $...,$ $x_{n},$ $y_{1},$
$y_{2},$ $...,$ $y_{n}$ be $2n$ points from $I.$ Assume that%
\begin{equation*}
\left\vert x_{1}-t\right\vert +\left\vert x_{2}-t\right\vert +...+\left\vert
x_{n}-t\right\vert \geq\left\vert y_{1}-t\right\vert +\left\vert
y_{2}-t\right\vert +...+\left\vert y_{n}-t\right\vert
\end{equation*}
holds for every $t\in\left\{
x_{1},x_{2},...,x_{n},y_{1},y_{2},...,y_{n}\right\} .$ Then,%
\begin{equation*}
f\left( x_{1}\right) +f\left( x_{2}\right) +...+f\left( x_{n}\right) \geq
f\left( y_{1}\right) +f\left( y_{2}\right) +...+f\left( y_{n}\right) .
\end{equation*}
\end{quote}

We will not need this result, but we will rather use its weighted version:

\begin{quote}
\textbf{Theorem 8b, the weighted Karamata inequality in symmetric form.} Let 
$f$ be a convex function from an interval $I\subseteq\mathbb{R}$ to $\mathbb{%
R},$ and let $N$ be a positive integer. Let $z_{1},$ $z_{2},$ $...,$ $z_{N}$
be $N$ points from $I,$ and let $w_{1},$ $w_{2},$ $...,$ $w_{N}$ be $N$
reals. Assume that%
\begin{equation}
\sum_{k=1}^{N}w_{k}=0,  \label{1}
\end{equation}
and that%
\begin{equation}
\sum_{k=1}^{N}w_{k}\left\vert z_{k}-t\right\vert \geq0\text{ holds for every 
}t\in\left\{ z_{1},z_{2},...,z_{N}\right\} .  \label{2}
\end{equation}
Then,%
\begin{equation}
\sum_{k=1}^{N}w_{k}f\left( z_{k}\right) \geq0.  \label{3}
\end{equation}
\end{quote}

It is very easy to conclude Theorem 8a from Theorem 8b by setting $N=2n$ and 
\begin{align*}
& z_{1}=x_{1},\ \ \ \ \ \ \ \ \ \ z_{2}=x_{2},\ \ \ \ \ \ \ \ \ \ ...,\ \ \
\ \ \ \ \ \ \ z_{n}=x_{n}; \\
& z_{n+1}=y_{1},\ \ \ \ \ \ \ \ \ \ z_{n+2}=y_{2},\ \ \ \ \ \ \ \ \ \ ...,\
\ \ \ \ \ \ \ \ \ z_{2n}=y_{n}; \\
& w_{1}=w_{2}=...=w_{n}=1;\ \ \ \ \ \ \ \ \ \ w_{n+1}=w_{n+2}=...=w_{2n}=-1,
\end{align*}
but as I said, we will never use Theorem 8a in this paper.

Time for \textit{a remark to readers familiar with majorization theory.} One
may wonder why I call the two results above "Karamata inequalities". In
fact, the Karamata inequality in its most known form claims:

\begin{quote}
\textbf{Theorem 9, the Karamata inequality.} Let $f$ be a convex function
from an interval $I\subseteq\mathbb{R}$ to $\mathbb{R},$ and let $n$ be a
positive integer. Let $x_{1},$ $x_{2},$ $...,$ $x_{n},$ $y_{1},$ $y_{2},$ $%
...,$ $y_{n}$ be $2n$ points from $I$ such that $\left(
x_{1},x_{2},...,x_{n}\right) \succ\left( y_{1},y_{2},...,y_{n}\right) .$
Then,%
\begin{equation*}
f\left( x_{1}\right) +f\left( x_{2}\right) +...+f\left( x_{n}\right) \geq
f\left( y_{1}\right) +f\left( y_{2}\right) +...+f\left( y_{n}\right) .
\end{equation*}
\end{quote}

According to [2], post \#11, Lemma 1, the condition $\left(
x_{1},x_{2},...,x_{n}\right) \succ\left( y_{1},y_{2},...,y_{n}\right) $
yields that $\left\vert x_{1}-t\right\vert +\left\vert x_{2}-t\right\vert
+...+\left\vert x_{n}-t\right\vert \geq\left\vert y_{1}-t\right\vert
+\left\vert y_{2}-t\right\vert +...+\left\vert y_{n}-t\right\vert $ holds
for every real $t$ - and thus, in particular, for every $t\in\left\{
z_{1},z_{2},...,z_{n}\right\} .$ Hence, whenever the condition of Theorem 9
holds, the condition of Theorem 8a holds as well. Thus, Theorem 9 follows
from Theorem 8a. With just a little more work, we could also derive Theorem
8a from Theorem 9, so that Theorems 8a and 9 are equivalent.

Note that Theorem 8b is more general than the Fuchs inequality (a more
well-known weighted version of the Karamata inequality). See [5] for a
generalization of majorization theory to weighted families of points
(apparently already known long time ago), with a different approach to this
fact.

As promised, here is a \textit{proof of Theorem 8b:} First, substituting $%
t=\max\left\{ z_{1},z_{2},...,z_{N}\right\} $ into (2) (it is clear that
this $t$ satisfies $t\in\left\{ z_{1},z_{2},...,z_{N}\right\} $), we get $%
z_{k}\leq t$ for every $k\in\left\{ 1,2,...,N\right\} ,$ so that $%
z_{k}-t\leq0$ and thus $\left\vert z_{k}-t\right\vert =-\left(
z_{k}-t\right) =t-z_{k}$ for every $k\in\left\{ 1,2,...,N\right\} ,$ and
thus (2) becomes $\sum\limits_{k=1}^{N}w_{k}\left( t-z_{k}\right) \geq0.$ In
other words, $t\sum\limits_{k=1}^{N}w_{k}-\sum\limits_{k=1}^{N}w_{k}z_{k}%
\geq0.$ In sight of $\sum\limits_{k=1}^{N}w_{k}=0,$ this rewrites as $%
t\cdot0-\sum\limits_{k=1}^{N}w_{k}z_{k}\geq0.$ Hence, $\sum%
\limits_{k=1}^{N}w_{k}z_{k}\leq0.$

Similarly, substituting $t=\min\left\{ z_{1},z_{2},...,z_{N}\right\} $ into
(2), we get $\sum\limits_{k=1}^{N}w_{k}z_{k}\geq0.$ Thus, $%
\sum\limits_{k=1}^{N}w_{k}z_{k}=0.$

The function $f:I\rightarrow\mathbb{R}$ is convex, and $z_{1},$ $z_{2},$ $%
...,$ $z_{N}$ are finitely many points from $I.$ Hence, Theorem 6 yields the
existence of two real constants $u$ and $v$ and $N$ nonnegative constants $%
a_{1},$ $a_{2},$ $...,$ $a_{N}$ such that%
\begin{equation*}
f\left( t\right) =vt+u+\sum\limits_{i=1}^{N}a_{i}\left\vert
t-z_{i}\right\vert \text{ holds for every }t\in\left\{
z_{1},z_{2},...,z_{N}\right\} .
\end{equation*}
Thus,%
\begin{equation*}
f\left( z_{k}\right) =vz_{k}+u+\sum\limits_{i=1}^{N}a_{i}\left\vert
z_{k}-z_{i}\right\vert \ \ \ \ \ \ \ \ \ \ \text{for every }k\in\left\{
1,2,...,N\right\} .
\end{equation*}
Hence,%
\begin{align*}
\sum_{k=1}^{N}w_{k}f\left( z_{k}\right) & =\sum_{k=1}^{N}w_{k}\left(
vz_{k}+u+\sum\limits_{i=1}^{N}a_{i}\left\vert z_{k}-z_{i}\right\vert \right)
=v\underbrace{\sum_{k=1}^{N}w_{k}z_{k}}_{=0}+u\underbrace{\sum_{k=1}^{N}w_{k}%
}_{=0}+\sum_{k=1}^{N}w_{k}\sum\limits_{i=1}^{N}a_{i}\left\vert
z_{k}-z_{i}\right\vert \\
& =\sum_{k=1}^{N}w_{k}\sum\limits_{i=1}^{N}a_{i}\left\vert
z_{k}-z_{i}\right\vert =\sum\limits_{i=1}^{N}a_{i}\underbrace{%
\sum_{k=1}^{N}w_{k}\left\vert z_{k}-z_{i}\right\vert }_{\geq0\text{
according to (2) for }t=z_{i}}\geq0.
\end{align*}
Thus, Theorem 8b is proven.

\begin{center}
\textbf{4. A property of zero-sum vectors}
\end{center}

Next, we are going to show some properties of real vectors.

If $k$ is an integer and $v\in\mathbb{R}^{k}$ is a vector, then, for any $%
i\in\left\{ 1,2,...,k\right\} ,$ we denote by $v_{i}$ the $i$-th coordinate
of the vector $v.$ Then, $v=\left( 
\begin{array}{c}
v_{1} \\ 
v_{2} \\ 
... \\ 
v_{k}%
\end{array}
\right) .$

Let $n$ be a positive integer. We consider the vector space $\mathbb{R}^{n}.$
Let $\left( e_{1},e_{2},...,e_{n}\right) $ be the standard basis of this
vector space $\mathbb{R}^{n};$ in other words, for every $i\in\left\{
1,2,...,n\right\} ,$ let $e_{i}$ be the vector from $\mathbb{R}^{n}$ such
that $\left( e_{i}\right) _{i}=1$ and $\left( e_{i}\right) _{j}=0$ for every 
$j\in\left\{ 1,2,...,n\right\} \setminus\left\{ i\right\} .$ Let $V_{n}$ be
the subspace of $\mathbb{R}^{n}$ defined by%
\begin{equation*}
V_{n}=\left\{ x\in\mathbb{R}^{n}\ \mid\ x_{1}+x_{2}+...+x_{n}=0\right\} .
\end{equation*}

For any $u\in\left\{ 1,2,...,n\right\} $ and any two \textit{distinct}
numbers $i$ and $j$ from the set $\left\{ 1,2,...,n\right\} ,$ we have%
\begin{equation}
\left( e_{i}-e_{j}\right) _{u}=\left\{ 
\begin{array}{c}
1,\text{ if }u=i; \\ 
-1,\text{ if }u=j; \\ 
0,\text{ if }u\neq i\text{ and }u\neq j%
\end{array}
\right. .  \label{4}
\end{equation}

Clearly, $e_{i}-e_{j}\in V_{n}$ for any two numbers $i$ and $j$ from the set 
$\left\{ 1,2,...,n\right\} .$

For any vector $t\in\mathbb{R}^{n},$ we denote $I\left( t\right) =\left\{
k\in\left\{ 1,2,...,n\right\} \mid t_{k}>0\right\} $ and $J\left( t\right)
=\left\{ k\in\left\{ 1,2,...,n\right\} \mid t_{k}<0\right\} .$ Obviously,
for every $t\in\mathbb{R}^{n},$ the sets $I\left( t\right) $ and $J\left(
t\right) $ are disjoint.

Now we are going to show:

\begin{quote}
\textbf{Theorem 10.} Let $n$ be a positive integer. Let $x\in V_{n}$ be a
vector. Then, there exist nonnegative reals $a_{i,j}$ for all pairs $\left(
i,j\right) \in I\left( x\right) \times J\left( x\right) $ such that%
\begin{equation*}
x=\sum_{\left( i,j\right) \in I\left( x\right) \times J\left( x\right)
}a_{i,j}\left( e_{i}-e_{j}\right) .
\end{equation*}
\end{quote}

\textit{Proof of Theorem 10.} We will prove Theorem 10 by induction over $%
\left\vert I\left( x\right) \right\vert +\left\vert J\left( x\right)
\right\vert .$

The \textit{basis of the induction} - the case when $\left\vert I\left(
x\right) \right\vert +\left\vert J\left( x\right) \right\vert =0$ - is
trivial: If $\left\vert I\left( x\right) \right\vert +\left\vert J\left(
x\right) \right\vert =0,$ then $I\left( x\right) =J\left( x\right)
=\varnothing$ and $x=0,$ so that $x=\sum\limits_{\left( i,j\right) \in
I\left( x\right) \times J\left( x\right) }a_{i,j}\left( e_{i}-e_{j}\right) $
holds because both sides of this equation are $0.$

Now we come to the \textit{induction step}: Let $r$ be a positive integer.
Assume that Theorem 10 holds for all $x\in V_{n}$ with $\left\vert I\left(
x\right) \right\vert +\left\vert J\left( x\right) \right\vert <r.$ We have
to show that Theorem 10 holds for all $x\in V_{n}$ with $\left\vert I\left(
x\right) \right\vert +\left\vert J\left( x\right) \right\vert =r.$

In order to prove this, we let $z\in V_{n}$ be an arbitrary vector with $%
\left\vert I\left( z\right) \right\vert +\left\vert J\left( z\right)
\right\vert =r.$ We then have to prove that Theorem 10 holds for $x=z.$ In
other words, we have to show that there exist nonnegative reals $a_{i,j}$
for all pairs $\left( i,j\right) \in I\left( z\right) \times J\left(
z\right) $ such that%
\begin{equation}
z=\sum_{\left( i,j\right) \in I\left( z\right) \times J\left( z\right)
}a_{i,j}\left( e_{i}-e_{j}\right) .  \label{5}
\end{equation}

First, $\left\vert I\left( z\right) \right\vert +\left\vert J\left( z\right)
\right\vert =r$ and $r>0$ yield $\left\vert I\left( z\right) \right\vert
+\left\vert J\left( z\right) \right\vert >0.$ Hence, at least one of the
sets $I\left( z\right) $ and $J\left( z\right) $ is non-empty.

Now, since $z\in V_{n},$ we have $z_{1}+z_{2}+...+z_{n}=0.$ Hence, either $%
z_{k}=0$ for every $k\in\left\{ 1,2,...,n\right\} ,$ or there is at least
one positive number and at least one negative number in the set $\left\{
z_{1},z_{2},...,z_{n}\right\} .$ The first case is impossible (since at
least one of the sets $I\left( z\right) $ and $J\left( z\right) $ is
non-empty). Thus, the second case must hold - i. e., there is at least one
positive number and at least one negative number in the set $\left\{
z_{1},z_{2},...,z_{n}\right\} .$ In other words, there exists a number $%
u\in\left\{ 1,2,...,n\right\} $ such that $z_{u}>0,$ and a number $%
v\in\left\{ 1,2,...,n\right\} $ such that $z_{v}<0.$ Of course, $z_{u}>0$
yields $u\in I\left( z\right) ,$ and $z_{v}<0$ yields $v\in J\left( z\right)
.$ Needless to say that $u\neq v.$

Now, we distinguish between two cases: the \textit{first case} will be the
case when $z_{u}+z_{v}\geq0,$ and the \textit{second case} will be the case
when $z_{u}+z_{v}\leq0.$

Let us consider the \textit{first case:} In this case, $z_{u}+z_{v}\geq0.$
Then, let $z^{\prime}=z+z_{v}\left( e_{u}-e_{v}\right) .$ Since $z\in V_{n}$
and $e_{u}-e_{v}\in V_{n},$ we have $z+z_{v}\left( e_{u}-e_{v}\right) \in
V_{n}$ (since $V_{n}$ is a vector space), so that $z^{\prime}\in V_{n}.$
From $z^{\prime}=z+z_{v}\left( e_{u}-e_{v}\right) ,$ the coordinate
representation of the vector $z^{\prime}$ is easily obtained:%
\begin{equation*}
z^{\prime}=\left( 
\begin{array}{c}
z_{1}^{\prime} \\ 
z_{2}^{\prime} \\ 
... \\ 
z_{n}^{\prime}%
\end{array}
\right) ,\ \ \ \ \ \ \ \ \ \ \text{where }\left\{ 
\begin{array}{c}
z_{k}^{\prime}=z_{k}\text{ for all }k\in\left\{ 1,2,...,n\right\}
\setminus\left\{ u,v\right\} ; \\ 
z_{u}^{\prime}=z_{u}+z_{v}; \\ 
z_{v}^{\prime}=0%
\end{array}
\right. .
\end{equation*}
It is readily seen from this that $I\left( z^{\prime}\right) \subseteq
I\left( z\right) $, so that $\left\vert I\left( z^{\prime}\right)
\right\vert \leq\left\vert I\left( z\right) \right\vert .$ Besides, $J\left(
z^{\prime}\right) \subseteq J\left( z\right) .$ Moreover, $J\left(
z^{\prime}\right) $ is a proper subset of $J\left( z\right) ,$ because $%
v\notin J\left( z^{\prime}\right) $ (since $z_{v}^{\prime}$ is not $<0,$ but 
$=0$) but $v\in J\left( z\right) .$ Hence, $\left\vert J\left(
z^{\prime}\right) \right\vert <\left\vert J\left( z\right) \right\vert .$
Combined with $\left\vert I\left( z^{\prime}\right) \right\vert
\leq\left\vert I\left( z\right) \right\vert ,$ this yields $\left\vert
I\left( z^{\prime}\right) \right\vert +\left\vert J\left( z^{\prime }\right)
\right\vert <\left\vert I\left( z\right) \right\vert +\left\vert J\left(
z\right) \right\vert .$ In view of $\left\vert I\left( z\right) \right\vert
+\left\vert J\left( z\right) \right\vert =r,$ this becomes $\left\vert
I\left( z^{\prime}\right) \right\vert +\left\vert J\left( z^{\prime}\right)
\right\vert <r.$ Thus, since we have assumed that Theorem 10 holds for all $%
x\in V_{n}$ with $\left\vert I\left( x\right) \right\vert +\left\vert
J\left( x\right) \right\vert <r,$ we can apply Theorem 10 to $x=z^{\prime},$
and we see that there exist nonnegative reals $a_{i,j}^{\prime}$ for all
pairs $\left( i,j\right) \in I\left( z^{\prime}\right) \times J\left(
z^{\prime}\right) $ such that%
\begin{equation*}
z^{\prime}=\sum_{\left( i,j\right) \in I\left( z^{\prime}\right) \times
J\left( z^{\prime}\right) }a_{i,j}^{\prime}\left( e_{i}-e_{j}\right) .
\end{equation*}
Now, $z^{\prime}=z+z_{v}\left( e_{u}-e_{v}\right) $ yields $z=z^{\prime
}-z_{v}\left( e_{u}-e_{v}\right) .$ Since $z_{v}<0,$ we have $-z_{v}>0,$ so
that, particularly, $-z_{v}$ is nonnegative.

Since $I\left( z^{\prime}\right) \subseteq I\left( z\right) $ and $J\left(
z^{\prime}\right) \subseteq J\left( z\right) ,$ we have $I\left(
z^{\prime}\right) \times J\left( z^{\prime}\right) \subseteq I\left(
z\right) \times J\left( z\right) .$ Also, $\left( u,v\right) \in I\left(
z\right) \times J\left( z\right) $ (because $u\in I\left( z\right) $ and $%
v\in J\left( z\right) $) and $\left( u,v\right) \notin I\left(
z^{\prime}\right) \times J\left( z^{\prime}\right) $ (because $v\notin
J\left( z^{\prime}\right) $).

Hence, the sets $I\left( z^{\prime}\right) \times J\left( z^{\prime }\right) 
$ and $\left\{ \left( u,v\right) \right\} $ are two disjoint subsets of the
set $I\left( z\right) \times J\left( z\right) .$ We can thus define
nonnegative reals $a_{i,j}$ for all pairs $\left( i,j\right) \in I\left(
z\right) \times J\left( z\right) $ by setting%
\begin{equation*}
a_{i,j}=\left\{ 
\begin{array}{c}
a_{i,j}^{\prime},\text{ if }\left( i,j\right) \in I\left( z^{\prime }\right)
\times J\left( z^{\prime}\right) ; \\ 
-z_{v},\text{ if }\left( i,j\right) =\left( u,v\right) ; \\ 
0,\text{ if neither of the two cases above holds}%
\end{array}
\right.
\end{equation*}
(these $a_{i,j}$ are all nonnegative because $a_{i,j}^{\prime},$ $-z_{v}$
and $0$ are nonnegative). Then,%
\begin{align*}
& \sum_{\left( i,j\right) \in I\left( z\right) \times J\left( z\right)
}a_{i,j}\left( e_{i}-e_{j}\right) \\
& =\sum_{\left( i,j\right) \in I\left( z^{\prime}\right) \times J\left(
z^{\prime}\right) }a_{i,j}^{\prime}\left( e_{i}-e_{j}\right) +\sum_{\left(
i,j\right) =\left( u,v\right) }\left( -z_{v}\right) \left(
e_{i}-e_{j}\right) +\left( \text{sum of some }0\text{'s}\right) \\
& =\sum_{\left( i,j\right) \in I\left( z^{\prime}\right) \times J\left(
z^{\prime}\right) }a_{i,j}^{\prime}\left( e_{i}-e_{j}\right) +\left(
-z_{v}\right) \left( e_{u}-e_{v}\right) +0=z^{\prime}+\left( -z_{v}\right)
\left( e_{u}-e_{v}\right) +0 \\
& =\left( z+z_{v}\left( e_{u}-e_{v}\right) \right) +\left( -z_{v}\right)
\left( e_{u}-e_{v}\right) +0=z.
\end{align*}
Thus, (5) is fulfilled.

Similarly, we can fulfill (5) in the \textit{second case}, repeating the
arguments we have done for the first case while occasionally interchanging $%
u $ with $v,$ as well as $I$ with $J,$ as well as $<$ with $>$. Here is a
brief outline of how we have to proceed in the second case: Denote $%
z^{\prime }=z-z_{u}\left( e_{u}-e_{v}\right) .$ Show that $z^{\prime}\in
V_{n}$ (as in the first case). Notice that%
\begin{equation*}
z^{\prime}=\left( 
\begin{array}{c}
z_{1}^{\prime} \\ 
z_{2}^{\prime} \\ 
... \\ 
z_{n}^{\prime}%
\end{array}
\right) ,\ \ \ \ \ \ \ \ \ \ \text{where }\left\{ 
\begin{array}{c}
z_{k}^{\prime}=z_{k}\text{ for all }k\in\left\{ 1,2,...,n\right\}
\setminus\left\{ u,v\right\} ; \\ 
z_{u}^{\prime}=0; \\ 
z_{v}^{\prime}=z_{u}+z_{v}%
\end{array}
\right. .
\end{equation*}
Prove that $u\notin I\left( z^{\prime}\right) $ (as we proved $v\notin
J\left( z^{\prime}\right) $ in the first case). Prove that $J\left(
z^{\prime}\right) \subseteq J\left( z\right) $ (similarly to the proof of $%
I\left( z^{\prime}\right) \subseteq I\left( z\right) $ in the first case)
and that $I\left( z^{\prime}\right) $ is a proper subset of $I\left(
z\right) $ (similarly to the proof that $J\left( z^{\prime}\right) $ is a
proper subset of $J\left( z\right) $ in the first case). Show that there
exist nonnegative reals $a_{i,j}^{\prime}$ for all pairs $\left( i,j\right)
\in I\left( z^{\prime}\right) \times J\left( z^{\prime}\right) $ such that%
\begin{equation*}
z^{\prime}=\sum_{\left( i,j\right) \in I\left( z^{\prime}\right) \times
J\left( z^{\prime}\right) }a_{i,j}^{\prime}\left( e_{i}-e_{j}\right)
\end{equation*}
(as in the first case). Note that $z_{u}$ is nonnegative (since $z_{u}>0$).
Prove that the sets $I\left( z^{\prime}\right) \times J\left( z^{\prime
}\right) $ and $\left\{ \left( u,v\right) \right\} $ are two disjoint
subsets of the set $I\left( z\right) \times J\left( z\right) $ (as in the
first case). Define nonnegative reals $a_{i,j}$ for all pairs $\left(
i,j\right) \in I\left( z\right) \times J\left( z\right) $ by setting 
\begin{equation*}
a_{i,j}=\left\{ 
\begin{array}{c}
a_{i,j}^{\prime},\text{ if }\left( i,j\right) \in I\left( z^{\prime }\right)
\times J\left( z^{\prime}\right) ; \\ 
z_{u},\text{ if }\left( i,j\right) =\left( u,v\right) ; \\ 
0,\text{ if neither of the two cases above holds}%
\end{array}
\right. .
\end{equation*}
Prove that these nonnegative reals $a_{i,j}$ fulfill (5).

Thus, in each of the two cases, we have proven that there exist nonnegative
reals $a_{i,j}$ for all pairs $\left( i,j\right) \in I\left( z\right) \times
J\left( z\right) $ such that (5) holds. Hence, Theorem 10 holds for $x=z.$
Thus, Theorem 10 is proven for all $x\in V_{n}$ with $\left\vert I\left(
x\right) \right\vert +\left\vert J\left( x\right) \right\vert =r.$ This
completes the induction step, and therefore, Theorem 10 is proven.

As an application of Theorem 10, we can now show:

\begin{quote}
\textbf{Theorem 11.} Let $n$ be a positive integer. Let $a_{1},$ $a_{2},$ $%
...,$ $a_{n}$ be $n$ nonnegative reals. Let $S$ be a finite set. For every $%
s\in S,$ let $r_{s}$ be an element of $\left( \mathbb{R}^{n}\right) ^{\ast }$
(in other words, a linear transformation from $\mathbb{R}^{n}$ to $\mathbb{R}
$), and let $b_{s}$ be a nonnegative real. Define a function $f:\mathbb{R}%
^{n}\rightarrow\mathbb{R}$ by%
\begin{equation*}
f\left( x\right) =\sum_{u=1}^{n}a_{u}\left\vert x_{u}\right\vert -\sum_{s\in
S}b_{s}\left\vert r_{s}x\right\vert ,\ \ \ \ \ \ \ \ \ \ \text{where }%
x=\left( 
\begin{array}{c}
x_{1} \\ 
x_{2} \\ 
... \\ 
x_{n}%
\end{array}
\right) \in\mathbb{R}^{n}.
\end{equation*}

Then, the following two assertions are equivalent:

\textit{Assertion }$\mathcal{A}_{1}$\textit{:} We have $f\left( x\right)
\geq0$ for every $x\in V_{n}.$

\textit{Assertion }$\mathcal{A}_{2}$\textit{:} We have $f\left(
e_{i}-e_{j}\right) \geq0$ for any two distinct integers $i$ and $j$ from $%
\left\{ 1,2,...,n\right\} .$
\end{quote}

\textit{Proof of Theorem 11.} We have to prove that the assertions $\mathcal{%
A}_{1}$ and $\mathcal{A}_{2}$ are equivalent.\ In other words, we have to
prove that $\mathcal{A}_{1}\Longrightarrow\mathcal{A}_{2}$ and $\mathcal{A}%
_{2}\Longrightarrow\mathcal{A}_{1}.$ Actually, $\mathcal{A}%
_{1}\Longrightarrow\mathcal{A}_{2}$ is trivial (we just have to use that $%
e_{i}-e_{j}\in V_{n}$ for any two numbers $i$ and $j$ from $\left\{
1,2,...,n\right\} $). It remains to show that $\mathcal{A}%
_{2}\Longrightarrow \mathcal{A}_{1}.$ So assume that Assertion $\mathcal{A}%
_{2}$ is valid, i. e. we have $f\left( e_{i}-e_{j}\right) \geq0$ for any two
distinct integers $i$ and $j$ from $\left\{ 1,2,...,n\right\} .$ We have to
prove that Assertion $\mathcal{A}_{1}$ holds, i. e. that $f\left( x\right)
\geq0$ for every $x\in V_{n}.$

So let $x\in V_{n}$ be some vector. According to Theorem 10, there exist
nonnegative reals $a_{i,j}$ for all pairs $\left( i,j\right) \in I\left(
x\right) \times J\left( x\right) $ such that%
\begin{equation*}
x=\sum_{\left( i,j\right) \in I\left( x\right) \times J\left( x\right)
}a_{i,j}\left( e_{i}-e_{j}\right) .
\end{equation*}
We will now show that%
\begin{equation}
\left\vert x_{u}\right\vert =\sum_{\left( i,j\right) \in I\left( x\right)
\times J\left( x\right) }a_{i,j}\left\vert \left( e_{i}-e_{j}\right)
_{u}\right\vert \ \ \ \ \ \ \ \ \ \ \text{for every }u\in\left\{
1,2,...,n\right\} .  \label{6}
\end{equation}
Here, of course, $\left( e_{i}-e_{j}\right) _{u}$ means the $u$-th
coordinate of the vector $e_{i}-e_{j}.$

In fact, two cases are possible: the case when $x_{u}\geq0,$ and the case
when $x_{u}<0.$ We will consider these cases separately.

\textit{Case 1:} We have $x_{u}\geq0.$ Then, $\left\vert x_{u}\right\vert
=x_{u}.$ Hence, in this case, we have $\left( e_{i}-e_{j}\right) _{u}\geq0$
for any two numbers $i\in I\left( x\right) $ and $j\in J\left( x\right) $
(in fact, $j\in J\left( x\right) $ yields $x_{j}<0,$ so that $u\neq j$
(because $x_{j}<0$ and $x_{u}\geq0$) and thus $\left( e_{j}\right) _{u}=0,$
so that $\left( e_{i}-e_{j}\right) _{u}=\left( e_{i}\right) _{u}-\left(
e_{j}\right) _{u}=\left( e_{i}\right) _{u}-0=\left( e_{i}\right)
_{u}=\left\{ 
\begin{array}{c}
1,\text{ if }u=i; \\ 
0,\text{ if }u\neq i%
\end{array}
\right. \geq0$). Thus, $\left( e_{i}-e_{j}\right) _{u}=\left\vert \left(
e_{i}-e_{j}\right) _{u}\right\vert $ for any two numbers $i\in I\left(
x\right) $ and $j\in J\left( x\right) .$ Thus, 
\begin{align*}
\left\vert x_{u}\right\vert & =x_{u}=\sum_{\left( i,j\right) \in I\left(
x\right) \times J\left( x\right) }a_{i,j}\left( e_{i}-e_{j}\right) _{u}\ \ \
\ \ \ \ \ \ \ \left( \text{since }x=\sum_{\left( i,j\right) \in I\left(
x\right) \times J\left( x\right) }a_{i,j}\left( e_{i}-e_{j}\right) \right) \\
& =\sum_{\left( i,j\right) \in I\left( x\right) \times J\left( x\right)
}a_{i,j}\left\vert \left( e_{i}-e_{j}\right) _{u}\right\vert ,
\end{align*}
and (6) is proven.

\textit{Case 2:} We have $x_{u}<0.$ Then, $u\in J\left( x\right) $ and $%
\left\vert x_{u}\right\vert =-x_{u}.$ Hence, in this case, we have $\left(
e_{i}-e_{j}\right) _{u}\leq0$ for any two numbers $i\in I\left( x\right) $
and $j\in J\left( x\right) $ (in fact, $i\in I\left( x\right) $ yields $%
x_{i}>0,$ so that $u\neq i$ (because $x_{i}>0$ and $x_{u}<0$) and thus $%
\left( e_{i}\right) _{u}=0,$ so that $\left( e_{i}-e_{j}\right) _{u}=\left(
e_{i}\right) _{u}-\left( e_{j}\right) _{u}=0-\left( e_{j}\right)
_{u}=-\left( e_{j}\right) _{u}=-\left\{ 
\begin{array}{c}
1,\text{ if }u=j; \\ 
0,\text{ if }u\neq j%
\end{array}
\right. \leq0$). Thus, $-\left( e_{i}-e_{j}\right) _{u}=\left\vert \left(
e_{i}-e_{j}\right) _{u}\right\vert $ for any two numbers $i\in I\left(
x\right) $ and $j\in J\left( x\right) .$ Thus, 
\begin{align*}
\left\vert x_{u}\right\vert & =-x_{u}=-\sum_{\left( i,j\right) \in I\left(
x\right) \times J\left( x\right) }a_{i,j}\left( e_{i}-e_{j}\right) _{u}\ \ \
\ \ \ \ \ \ \ \left( \text{since }x=\sum_{\left( i,j\right) \in I\left(
x\right) \times J\left( x\right) }a_{i,j}\left( e_{i}-e_{j}\right) \right) \\
& =\sum_{\left( i,j\right) \in I\left( x\right) \times J\left( x\right)
}a_{i,j}\left( -\left( e_{i}-e_{j}\right) _{u}\right) =\sum_{\left(
i,j\right) \in I\left( x\right) \times J\left( x\right) }a_{i,j}\left\vert
\left( e_{i}-e_{j}\right) _{u}\right\vert ,
\end{align*}
and (6) is proven.

Hence, in both cases, (6) is proven. Thus, (6) always holds. Now let us
continue our proof of $\mathcal{A}_{2}\Longrightarrow\mathcal{A}_{1}$:

We have%
\begin{align*}
\sum_{s\in S}b_{s}\left\vert r_{s}x\right\vert & =\sum_{s\in
S}b_{s}\left\vert r_{s}\sum_{\left( i,j\right) \in I\left( x\right) \times
J\left( x\right) }a_{i,j}\left( e_{i}-e_{j}\right) \right\vert \ \ \ \ \ \ \
\ \ \ \left( \text{since }x=\sum_{\left( i,j\right) \in I\left( x\right)
\times J\left( x\right) }a_{i,j}\left( e_{i}-e_{j}\right) \right) \\
& =\sum_{s\in S}b_{s}\left\vert \sum_{\left( i,j\right) \in I\left( x\right)
\times J\left( x\right) }a_{i,j}r_{s}\left( e_{i}-e_{j}\right) \right\vert \\
& \leq\sum_{s\in S}b_{s}\sum_{\left( i,j\right) \in I\left( x\right) \times
J\left( x\right) }a_{i,j}\left\vert r_{s}\left( e_{i}-e_{j}\right)
\right\vert \\
& \ \ \ \ \ \ \ \ \ \ \left( \text{by the triangle inequality, since all }%
a_{i,j}\text{ and all }b_{s}\text{ are nonnegative}\right) .
\end{align*}
Thus,%
\begin{align*}
f\left( x\right) & =\sum_{u=1}^{n}a_{u}\left\vert x_{u}\right\vert
-\sum_{s\in S}b_{s}\left\vert r_{s}x\right\vert
\geq\sum_{u=1}^{n}a_{u}\left\vert x_{u}\right\vert -\sum_{s\in
S}b_{s}\sum_{\left( i,j\right) \in I\left( x\right) \times J\left( x\right)
}a_{i,j}\left\vert r_{s}\left( e_{i}-e_{j}\right) \right\vert \\
& =\sum_{u=1}^{n}a_{u}\cdot\sum_{\left( i,j\right) \in I\left( x\right)
\times J\left( x\right) }a_{i,j}\left\vert \left( e_{i}-e_{j}\right)
_{u}\right\vert -\sum_{s\in S}b_{s}\sum_{\left( i,j\right) \in I\left(
x\right) \times J\left( x\right) }a_{i,j}\left\vert r_{s}\left(
e_{i}-e_{j}\right) \right\vert \ \ \ \ \ \ \ \ \ \ \left( \text{by (6)}%
\right) \\
& =\sum_{\left( i,j\right) \in I\left( x\right) \times J\left( x\right)
}a_{i,j}\sum_{u=1}^{n}a_{u}\left\vert \left( e_{i}-e_{j}\right)
_{u}\right\vert -\sum_{\left( i,j\right) \in I\left( x\right) \times J\left(
x\right) }a_{i,j}\sum_{s\in S}b_{s}\left\vert r_{s}\left( e_{i}-e_{j}\right)
\right\vert \\
& =\sum_{\left( i,j\right) \in I\left( x\right) \times J\left( x\right)
}a_{i,j}\cdot\left( \sum_{u=1}^{n}a_{u}\left\vert \left( e_{i}-e_{j}\right)
_{u}\right\vert -\sum_{s\in S}b_{s}\left\vert r_{s}\left( e_{i}-e_{j}\right)
\right\vert \right) \\
& =\sum_{\left( i,j\right) \in I\left( x\right) \times J\left( x\right) }%
\underbrace{a_{i,j}}_{\geq0}\cdot\underbrace{f\left( e_{i}-e_{j}\right) }%
_{\geq0}\geq0.
\end{align*}
(Here, $f\left( e_{i}-e_{j}\right) \geq0$ because $i$ and $j$ are two
distinct integers from $\left\{ 1,2,...,n\right\} ;$ in fact, $i$ and $j$
are distinct because $i\in I\left( x\right) $ and $j\in J\left( x\right) ,$
and the sets $I\left( x\right) $ and $J\left( x\right) $ are disjoint.)

Hence, we have obtained $f\left( x\right) \geq0.$ This proves the assertion $%
\mathcal{A}_{1}.$ Therefore, the implication $\mathcal{A}_{2}\Longrightarrow 
\mathcal{A}_{1}$ is proven, and the proof of Theorem 11 is complete.

\begin{center}
\textbf{5. Restating Theorem 11}
\end{center}

Now we consider a result which follows from Theorem 11 pretty obviously
(although the formalization of the proof is going to be gruelling):

\begin{quote}
\textbf{Theorem 12.} Let $n$ be a nonnegative integer. Let $a_{1},$ $a_{2},$ 
$...,$ $a_{n}$ and $a$ be $n+1$ nonnegative reals. Let $S$ be a finite set.
For every $s\in S,$ let $r_{s}$ be an element of $\left( \mathbb{R}%
^{n}\right) ^{\ast}$ (in other words, a linear transformation from $\mathbb{R%
}^{n}$ to $\mathbb{R}$), and let $b_{s}$ be a nonnegative real. Define a
function $g:\mathbb{R}^{n}\rightarrow\mathbb{R}$ by%
\begin{equation*}
g\left( x\right) =\sum_{u=1}^{n}a_{u}\left\vert x_{u}\right\vert
+a\left\vert x_{1}+x_{2}+...+x_{n}\right\vert -\sum_{s\in S}b_{s}\left\vert
r_{s}x\right\vert ,\ \ \ \ \ \ \ \ \ \ \text{where }x=\left( 
\begin{array}{c}
x_{1} \\ 
x_{2} \\ 
... \\ 
x_{n}%
\end{array}
\right) \in\mathbb{R}^{n}.
\end{equation*}

Then, the following two assertions are equivalent:

\textit{Assertion }$\mathcal{B}_{1}$\textit{:} We have $g\left( x\right)
\geq0$ for every $x\in\mathbb{R}^{n}.$

\textit{Assertion }$\mathcal{B}_{2}$\textit{:} We have $g\left( e_{i}\right)
\geq0$ for every integer $i\in\left\{ 1,2,...,n\right\} ,$ and $g\left(
e_{i}-e_{j}\right) \geq0$ for any two distinct integers $i$ and $j$ from $%
\left\{ 1,2,...,n\right\} .$
\end{quote}

\textit{Proof of Theorem 12.} We are going to restate Theorem 12 before we
actually prove it. But first, we introduce a notation:

Let $\left( \widetilde{e_{1}},\widetilde{e_{2}},...,\widetilde{e_{n-1}}%
\right) $ be the standard basis of the vector space $\mathbb{R}^{n-1};$ in
other words, for every $i\in\left\{ 1,2,...,n-1\right\} ,$ let $\widetilde{%
e_{i}}$ be the vector from $\mathbb{R}^{n-1}$ such that $\left( \widetilde{%
e_{i}}\right) _{i}=1$ and $\left( \widetilde{e_{i}}\right) _{j}=0$ for every 
$j\in\left\{ 1,2,...,n-1\right\} \setminus\left\{ i\right\} .$

Now we will restate Theorem 12 by renaming $n$ into $n-1$ (thus replacing $%
e_{i}$ by $\widetilde{e_{i}}$ as well) and $a$ into $a_{n}$:

\begin{quote}
\textbf{Theorem 12b.} Let $n$ be a positive integer. Let $a_{1},$ $a_{2},$ $%
...,$ $a_{n-1},$ $a_{n}$ be $n$ nonnegative reals. Let $S$ be a finite set.
For every $s\in S,$ let $r_{s}$ be an element of $\left( \mathbb{R}%
^{n-1}\right) ^{\ast}$ (in other words, a linear transformation from $%
\mathbb{R}^{n-1}$ to $\mathbb{R}$), and let $b_{s}$ be a nonnegative real.
Define a function $g:\mathbb{R}^{n-1}\rightarrow\mathbb{R}$ by%
\begin{equation*}
g\left( x\right) =\sum_{u=1}^{n-1}a_{u}\left\vert x_{u}\right\vert
+a_{n}\left\vert x_{1}+x_{2}+...+x_{n-1}\right\vert -\sum_{s\in
S}b_{s}\left\vert r_{s}x\right\vert ,\ \ \ \ \ \ \ \ \ \ \text{where }%
x=\left( 
\begin{array}{c}
x_{1} \\ 
x_{2} \\ 
... \\ 
x_{n-1}%
\end{array}
\right) \in\mathbb{R}^{n-1}.
\end{equation*}

Then, the following two assertions are equivalent:

\textit{Assertion }$\mathcal{C}_{1}$\textit{:} We have $g\left( x\right)
\geq0$ for every $x\in\mathbb{R}^{n-1}.$

\textit{Assertion }$\mathcal{C}_{2}$\textit{:} We have $g\left( \widetilde{%
e_{i}}\right) \geq0$ for every integer $i\in\left\{ 1,2,...,n-1\right\} ,$
and $g\left( \widetilde{e_{i}}-\widetilde{e_{j}}\right) \geq0$ for any two
distinct integers $i$ and $j$ from $\left\{ 1,2,...,n-1\right\} .$
\end{quote}

Theorem 12b is equivalent to Theorem 12 (because Theorem 12b is just Theorem
12, applied to $n-1$ instead of $n$). Thus, proving Theorem 12b will be
enough to verify Theorem 12.

\textit{Proof of Theorem 12b.} The implication $\mathcal{C}%
_{1}\Longrightarrow \mathcal{C}_{2}$ is absolutely trivial. Hence, in order
to establish Theorem 12b, it only remains to prove the implication $\mathcal{%
C}_{2}\Longrightarrow \mathcal{C}_{1}.$

So assume that the assertion $\mathcal{C}_{2}$ holds, i. e. that we have $%
g\left( \widetilde{e_{i}}\right) \geq0$ for every integer $i\in\left\{
1,2,...,n-1\right\} ,$ and $g\left( \widetilde{e_{i}}-\widetilde{e_{j}}%
\right) \geq0$ for any two distinct integers $i$ and $j$ from $\left\{
1,2,...,n-1\right\} .$ We want to show that Assertion $\mathcal{C}_{1}$
holds, i. e. that $g\left( x\right) \geq0$ is satisfied for every $x\in%
\mathbb{R}^{n-1}.$

Since $\left( \widetilde{e_{1}},\widetilde{e_{2}},...,\widetilde{e_{n-1}}%
\right) $ is the standard basis of the vector space $\mathbb{R}^{n-1},$
every vector $x\in\mathbb{R}^{n-1}$ satisfies $x=\sum\limits_{i=1}^{n-1}x_{i}%
\widetilde{e_{i}}.$

Since $\left( e_{1},e_{2},...,e_{n}\right) $ is the standard basis of the
vector space $\mathbb{R}^{n},$ every vector $x\in\mathbb{R}^{n}$ satisfies $%
x=\sum\limits_{i=1}^{n}x_{i}e_{i}.$

Let $\phi_{n}:\mathbb{R}^{n-1}\rightarrow\mathbb{R}^{n}$ be the linear
transformation defined by $\phi_{n}\widetilde{e_{i}}=e_{i}-e_{n}$ for every $%
i\in\left\{ 1,2,...,n-1\right\} .$ (This linear transformation is uniquely
defined this way because $\left( \widetilde{e_{1}},\widetilde{e_{2}},...,%
\widetilde{e_{n-1}}\right) $ is a basis of $\mathbb{R}^{n-1}.$) For every $%
x\in\mathbb{R}^{n-1},$ we then have%
\begin{align}
\phi_{n}x & =\phi_{n}\left( \sum\limits_{i=1}^{n-1}x_{i}\widetilde{e_{i}}%
\right) =\sum\limits_{i=1}^{n-1}x_{i}\phi_{n}\widetilde{e_{i}}\ \ \ \ \ \ \
\ \ \ \left( \text{since }\phi_{n}\text{ is linear}\right)  \notag \\
& =\sum\limits_{i=1}^{n-1}x_{i}\left( e_{i}-e_{n}\right) =\sum
\limits_{i=1}^{n-1}x_{i}e_{i}-\sum\limits_{i=1}^{n-1}x_{i}e_{n}=\sum
\limits_{i=1}^{n-1}x_{i}e_{i}-\left( x_{1}+x_{2}+...+x_{n-1}\right) e_{n} 
\notag \\
& =\left( 
\begin{array}{c}
x_{1} \\ 
x_{2} \\ 
... \\ 
x_{n-1} \\ 
-\left( x_{1}+x_{2}+...+x_{n-1}\right)%
\end{array}
\right) ,  \label{7}
\end{align}
Consequently, $\phi_{n}x\in V_{n}$ for every $x\in\mathbb{R}^{n-1}$. Hence, $%
\func{Im}\phi_{n}\subseteq V_{n}.$

Let $\psi_{n}:\mathbb{R}^{n}\rightarrow\mathbb{R}^{n-1}$ be the linear
transformation defined by $\psi_{n}e_{i}=\left\{ 
\begin{array}{c}
\widetilde{e_{i}},\text{ if }i\in\left\{ 1,2,...,n-1\right\} ; \\ 
0,\text{ if }i=n%
\end{array}
\right. $ for every $i\in\left\{ 1,2,...,n\right\} .$ (This linear
transformation is uniquely defined this way because $\left(
e_{1},e_{2},...,e_{n}\right) $ is a basis of $\mathbb{R}^{n}.$) For every $%
x\in\mathbb{R}^{n},$ we then have%
\begin{align*}
\psi_{n}x & =\psi_{n}\left( \sum\limits_{i=1}^{n}x_{i}e_{i}\right)
=\sum\limits_{i=1}^{n}x_{i}\psi_{n}e_{i}\ \ \ \ \ \ \ \ \ \ \left( \text{%
since }\psi_{n}\text{ is linear}\right) \\
& =\sum\limits_{i=1}^{n}x_{i}\left\{ 
\begin{array}{c}
\widetilde{e_{i}},\text{ if }i\in\left\{ 1,2,...,n-1\right\} ; \\ 
0,\text{ if }i=n%
\end{array}
\right. \\
& =\sum_{i=1}^{n-1}x_{i}\widetilde{e_{i}}=\left( 
\begin{array}{c}
x_{1} \\ 
x_{2} \\ 
... \\ 
x_{n-1}%
\end{array}
\right) .
\end{align*}

Then, $\psi_{n}\phi_{n}=\limfunc{id}$ (in fact, for every $i\in\left\{
1,2,...,n-1\right\} ,$ we have%
\begin{align*}
\psi_{n}\phi_{n}\widetilde{e_{i}} & =\psi_{n}\left( e_{i}-e_{n}\right)
=\psi_{n}e_{i}-\psi_{n}e_{n}\ \ \ \ \ \ \ \ \ \ \left( \text{since }\psi _{n}%
\text{ is linear}\right) \\
& =\widetilde{e_{i}}-0=\widetilde{e_{i}};
\end{align*}
thus, for every $x\in\mathbb{R}^{n-1},$ we have%
\begin{align*}
\psi_{n}\phi_{n}x & =\psi_{n}\phi_{n}\left( \sum\limits_{i=1}^{n-1}x_{i}%
\widetilde{e_{i}}\right) =\sum\limits_{i=1}^{n-1}x_{i}\psi_{n}\phi _{n}%
\widetilde{e_{i}} \\
& \ \ \ \ \ \ \ \ \ \ \left( \text{since the function }\psi_{n}\phi _{n}%
\text{ is linear, because }\psi_{n}\text{ and }\phi_{n}\text{ are linear}%
\right) \\
& =\sum\limits_{i=1}^{n-1}x_{i}\widetilde{e_{i}}=x,
\end{align*}
and therefore $\psi_{n}\phi_{n}=\limfunc{id}$).

We define a function $f:\mathbb{R}^{n}\rightarrow\mathbb{R}$ by%
\begin{equation*}
f\left( x\right) =\sum_{u=1}^{n}a_{u}\left\vert x_{u}\right\vert -\sum_{s\in
S}b_{s}\left\vert r_{s}\psi_{n}x\right\vert ,\ \ \ \ \ \ \ \ \ \ \text{where 
}x=\left( 
\begin{array}{c}
x_{1} \\ 
x_{2} \\ 
... \\ 
x_{n}%
\end{array}
\right) \in\mathbb{R}^{n}.
\end{equation*}
Note that%
\begin{equation}
f\left( -x\right) =f\left( x\right) \ \ \ \ \ \ \ \ \ \ \text{for every }x\in%
\mathbb{R}^{n},  \label{8}
\end{equation}
since%
\begin{align*}
f\left( -x\right) & =\sum_{u=1}^{n}a_{u}\left\vert \left( -x\right)
_{u}\right\vert -\sum_{s\in S}b_{s}\left\vert r_{s}\psi_{n}\left( -x\right)
\right\vert =\sum_{u=1}^{n}a_{u}\left\vert -x_{u}\right\vert -\sum_{s\in
S}b_{s}\left\vert -r_{s}\psi_{n}x\right\vert \\
& \ \ \ \ \ \ \ \ \ \ \left( \text{here, we have }r_{s}\psi_{n}\left(
-x\right) =-r_{s}\psi_{n}x\text{ since }r_{s}\text{ and }\psi_{n}\text{ are
linear functions}\right) \\
& =\sum_{u=1}^{n}a_{u}\left\vert x_{u}\right\vert -\sum_{s\in
S}b_{s}\left\vert r_{s}\psi_{n}x\right\vert =f\left( x\right) .
\end{align*}
Furthermore, I claim that%
\begin{equation}
f\left( \phi_{n}x\right) =g\left( x\right) \ \ \ \ \ \ \ \ \ \ \text{for
every }x\in\mathbb{R}^{n-1}.  \label{9}
\end{equation}
In order to prove this, we note that (7) yields $\left( \phi_{n}x\right)
_{u}=x_{u}$ for all $u\in\left\{ 1,2,...,n-1\right\} $ and $\left( \phi
_{n}x\right) _{n}=-\left( x_{1}+x_{2}+...+x_{n-1}\right) ,$ while $\psi
_{n}\phi_{n}=\limfunc{id}$ yields $\psi_{n}\phi_{n}x=x,$ so that%
\begin{align*}
f\left( \phi_{n}x\right) & =\sum_{u=1}^{n}a_{u}\left\vert \left( \phi
_{n}x\right) _{u}\right\vert -\sum_{s\in S}b_{s}\left\vert
r_{s}\psi_{n}\phi_{n}x\right\vert \\
& =\sum_{u=1}^{n-1}a_{u}\left\vert \left( \phi_{n}x\right) _{u}\right\vert
+a_{n}\left\vert \left( \phi_{n}x\right) _{n}\right\vert -\sum_{s\in
S}b_{s}\left\vert r_{s}\psi_{n}\phi_{n}x\right\vert \\
& =\sum_{u=1}^{n-1}a_{u}\left\vert x_{u}\right\vert +a_{n}\left\vert -\left(
x_{1}+x_{2}+...+x_{n-1}\right) \right\vert -\sum_{s\in S}b_{s}\left\vert
r_{s}x\right\vert \\
& \ \ \ \ \ \ \ \ \ \ \left( \text{since }\left( \phi_{n}x\right) _{u}=x_{u}%
\text{ for all }u\in\left\{ 1,2,...,n-1\right\} \text{ and}\right. \\
& \ \ \ \ \ \ \ \ \ \ \left. \left( \phi_{n}x\right) _{n}=-\left(
x_{1}+x_{2}+...+x_{n-1}\right) ,\text{ and }\psi_{n}\phi_{n}x=x\right) \\
& =\sum_{u=1}^{n-1}a_{u}\left\vert x_{u}\right\vert +a_{n}\left\vert
x_{1}+x_{2}+...+x_{n-1}\right\vert -\sum_{s\in S}b_{s}\left\vert
r_{s}x\right\vert =g\left( x\right) ,
\end{align*}
and thus (9) is proven.

Now, we are going to show that%
\begin{equation}
f\left( e_{i}-e_{j}\right) \geq0\text{ for any two distinct integers }i\text{
and }j\text{ from }\left\{ 1,2,...,n\right\} .  \label{10}
\end{equation}
In order to prove (10), we distinguish between three different cases:

\textit{Case 1:} We have $i\in\left\{ 1,2,...,n-1\right\} $ and $j\in\left\{
1,2,...,n-1\right\} .$

\textit{Case 2:} We have $i\in\left\{ 1,2,...,n-1\right\} $ and $j=n.$

\textit{Case 3:} We have $i=n$ and $j\in\left\{ 1,2,...,n-1\right\} .$

(In fact, the case when both $i=n$ and $j=n$ cannot occur, since $i$ and $j$
must be distinct).

In Case 1, we have%
\begin{align*}
f\left( e_{i}-e_{j}\right) & =f\left( \left( e_{i}-e_{n}\right) -\left(
e_{j}-e_{n}\right) \right) =f\left( \phi_{n}\widetilde{e_{i}}-\phi_{n}%
\widetilde{e_{j}}\right) \\
& =f\left( \phi_{n}\left( \widetilde{e_{i}}-\widetilde{e_{j}}\right) \right)
\ \ \ \ \ \ \ \ \ \ \left( \text{since }\phi_{n}\widetilde{e_{i}}-\phi_{n}%
\widetilde{e_{j}}=\phi_{n}\left( \widetilde{e_{i}}-\widetilde {e_{j}}\right)
,\text{ because }\phi_{n}\text{ is linear}\right) \\
& =g\left( \widetilde{e_{i}}-\widetilde{e_{j}}\right) \ \ \ \ \ \ \ \ \ \
\left( \text{after (9)}\right) \\
& \geq0\ \ \ \ \ \ \ \ \ \ \left( \text{by assumption}\right) .
\end{align*}

In Case 2, we have%
\begin{align*}
f\left( e_{i}-e_{j}\right) & =f\left( e_{i}-e_{n}\right) =f\left( \phi_{n}%
\widetilde{e_{i}}\right) =g\left( \widetilde{e_{i}}\right) \ \ \ \ \ \ \ \ \
\ \left( \text{after (9)}\right) \\
& \geq0\ \ \ \ \ \ \ \ \ \ \left( \text{by assumption}\right) .
\end{align*}

In Case 3, we have%
\begin{align*}
f\left( e_{i}-e_{j}\right) & =f\left( e_{n}-e_{j}\right) =f\left( -\left(
e_{j}-e_{n}\right) \right) =f\left( e_{j}-e_{n}\right) \ \ \ \ \ \ \ \ \ \
\left( \text{after (8)}\right) \\
& =f\left( \phi_{n}\widetilde{e_{j}}\right) =g\left( \widetilde{e_{j}}%
\right) \ \ \ \ \ \ \ \ \ \ \left( \text{after (9)}\right) \\
& \geq0\ \ \ \ \ \ \ \ \ \ \left( \text{by assumption}\right) .
\end{align*}

Thus, $f\left( e_{i}-e_{j}\right) \geq0$ holds in all three possible cases.
Hence, (10) is proven.

Now, our function $f:\mathbb{R}^{n}\rightarrow\mathbb{R}$ is defined by%
\begin{equation*}
f\left( x\right) =\sum_{u=1}^{n}a_{u}\left\vert x_{u}\right\vert -\sum_{s\in
S}b_{s}\left\vert r_{s}\psi_{n}x\right\vert ,\ \ \ \ \ \ \ \ \ \ \text{where 
}x=\left( 
\begin{array}{c}
x_{1} \\ 
x_{2} \\ 
... \\ 
x_{n}%
\end{array}
\right) \in\mathbb{R}^{n}.
\end{equation*}
Here, $n$ is a positive integer; the numbers $a_{1},$ $a_{2},$ $...,$ $a_{n}$
are $n$ nonnegative reals; the set $S$ is a finite set; for every $s\in S,$
the function $r_{s}\psi_{n}$ is an element of $\left( \mathbb{R}^{n}\right)
^{\ast}$ (in other words, a linear transformation from $\mathbb{R}^{n}$ to $%
\mathbb{R}$), and $b_{s}$ is a nonnegative real.

Hence, we can apply Theorem 11 to our function $f,$ and we obtain that for
our function $f,$ the Assertions $\mathcal{A}_{1}$ and $\mathcal{A}_{2}$ are
equivalent. In other words, our function $f$ satisfies Assertion $\mathcal{A}%
_{1}$ if and only if it satisfies Assertion $\mathcal{A}_{2}.$

Now, according to (10), our function $f$ satisfies Assertion $\mathcal{A}%
_{2}.$ Thus, this function $f$ must also satisfy Assertion $\mathcal{A}_{1}.$
In other words, $f\left( x\right) \geq0$ holds for every $x\in V_{n}.$
Hence, $f\left( \phi_{n}x\right) \geq0$ holds for every $x\in\mathbb{R}%
^{n-1} $ (because $\phi_{n}x\in V_{n},$ since $\func{Im}\phi _{n}\subseteq
V_{n}$). Since $f\left( \phi_{n}x\right) =g\left( x\right) $ according to
(9), we have therefore proven that $g\left( x\right) \geq0$ holds for every $%
x\in\mathbb{R}^{n-1}.$ Hence, Assertion $\mathcal{C}_{1}$ is proven. Thus,
we have showed that $\mathcal{C}_{2}\Longrightarrow \mathcal{C}_{1},$ and
thus the proof of Theorem 12b is complete.

Since Theorem 12b is equivalent to Theorem 12, this also proves Theorem 12.

As if this wasn't enough, here comes a further restatement of Theorem 12:

\begin{quote}
\textbf{Theorem 13.} Let $n$ be a nonnegative integer. Let $a_{1},$ $a_{2},$ 
$...,$ $a_{n}$ and $a$ be $n+1$ nonnegative reals. Let $S$ be a finite set.
For every $s\in S,$ let $r_{s,1},$ $r_{s,2},$ $...,$ $r_{s,n}$ be $n$
nonnegative reals, and let $b_{s}$ be a nonnegative real. Assume that the
following two assertions hold:%
\begin{align*}
a_{i}+a & \geq\sum_{s\in S}b_{s}r_{s,i}\ \ \ \ \ \ \ \ \ \ \text{for every }%
i\in\left\{ 1,2,...,n\right\} ; \\
a_{i}+a_{j} & \geq\sum_{s\in S}b_{s}\left\vert r_{s,i}-r_{s,j}\right\vert \
\ \ \ \ \ \ \ \ \ \text{for any two distinct integers }i\text{ and }j\text{
from }\left\{ 1,2,...,n\right\} .
\end{align*}
Let $y_{1},$ $y_{2},$ $...,$ $y_{n}$ be $n$ reals. Then,%
\begin{equation*}
\sum\limits_{i=1}^{n}a_{i}\left\vert y_{i}\right\vert +a\left\vert
\sum\limits_{v=1}^{n}y_{v}\right\vert -\sum\limits_{s\in S}b_{s}\left\vert
\sum\limits_{v=1}^{n}r_{s,v}y_{v}\right\vert \geq0.
\end{equation*}
\end{quote}

\textit{Proof of Theorem 13.} For every $s\in S,$ let $r_{s}=\left(
r_{s,1},r_{s,2},...,r_{s,n}\right) \in\left( \mathbb{R}^{n}\right) ^{\ast}$
be the $n$-dimensional covector whose $i$-th coordinate is $r_{s,i}$ for
every $i\in\left\{ 1,2,...,n\right\} .$ Define a function $g:\mathbb{R}%
^{n}\rightarrow\mathbb{R}$ by%
\begin{equation*}
g\left( x\right) =\sum_{u=1}^{n}a_{u}\left\vert x_{u}\right\vert
+a\left\vert x_{1}+x_{2}+...+x_{n}\right\vert -\sum_{s\in S}b_{s}\left\vert
r_{s}x\right\vert ,\ \ \ \ \ \ \ \ \ \ \text{where }x=\left( 
\begin{array}{c}
x_{1} \\ 
x_{2} \\ 
... \\ 
x_{n}%
\end{array}
\right) \in\mathbb{R}^{n}.
\end{equation*}

For every $i\in\left\{ 1,2,...,n\right\} ,$ we have $\left( e_{i}\right)
_{u}=\left\{ 
\begin{array}{c}
1,\text{ if }u=i; \\ 
0,\text{ if }u\neq i%
\end{array}
\right. $ for all $u\in\left\{ 1,2,...,n\right\} ,$ so that $\left(
e_{i}\right) _{1}+\left( e_{i}\right) _{2}+...+\left( e_{i}\right) _{n}=1,$
and for every $s\in S,$ we have%
\begin{align*}
r_{s}e_{i} & =\sum_{u=1}^{n}r_{s,u}\left( e_{i}\right) _{u}\ \ \ \ \ \ \ \ \
\ \left( \text{since }r_{s}=\left( r_{s,1},r_{s,2},...,r_{s,n}\right) \right)
\\
& =\sum_{u=1}^{n}r_{s,u}\left\{ 
\begin{array}{c}
1,\text{ if }u=i; \\ 
0,\text{ if }u\neq i%
\end{array}
\right. =r_{s,i},
\end{align*}
so that%
\begin{align*}
g\left( e_{i}\right) & =\sum_{u=1}^{n}a_{u}\left\vert \left( e_{i}\right)
_{u}\right\vert +a\left\vert \left( e_{i}\right) _{1}+\left( e_{i}\right)
_{2}+...+\left( e_{i}\right) _{n}\right\vert -\sum_{s\in S}b_{s}\left\vert
r_{s}e_{i}\right\vert \\
& =\sum_{u=1}^{n}a_{u}\left\vert \left\{ 
\begin{array}{c}
1,\text{ if }u=i; \\ 
0,\text{ if }u\neq i%
\end{array}
\right. \right\vert +a\left\vert 1\right\vert -\sum_{s\in S}b_{s}\left\vert
r_{s,i}\right\vert \\
& =a_{i}\left\vert 1\right\vert +a\left\vert 1\right\vert -\sum_{s\in S}b_{s}%
\underbrace{\left\vert r_{s,i}\right\vert }_{\substack{ =r_{s,i},  \\ \text{%
since}  \\ r_{s,i}\geq0}}=a_{i}+a-\sum_{s\in S}b_{s}r_{s,i}\geq0
\end{align*}
(since $a_{i}+a\geq\sum\limits_{s\in S}b_{s}r_{s,i}$ by the conditions of
Theorem 13).

For any two distinct integers $i$ and $j$ from $\left\{ 1,2,...,n\right\} ,$
we have $\left( e_{i}-e_{j}\right) _{u}=\left\{ 
\begin{array}{c}
1,\text{ if }u=i; \\ 
-1,\text{ if }u=j; \\ 
0,\text{ if }u\neq i\text{ and }u\neq j%
\end{array}
\right. \ $for all $u\in\left\{ 1,2,...,n\right\} ,$ so that $\left(
e_{i}-e_{j}\right) _{1}+\left( e_{i}-e_{j}\right) _{2}+...+\left(
e_{i}-e_{j}\right) _{n}=0,$ and for every $s\in S,$ we have%
\begin{align*}
r_{s}\left( e_{i}-e_{j}\right) & =\sum_{u=1}^{n}r_{s,u}\left(
e_{i}-e_{j}\right) _{u}\ \ \ \ \ \ \ \ \ \ \left( \text{since }r_{s}=\left(
r_{s,1},r_{s,2},...,r_{s,n}\right) \right) \\
& =\sum_{u=1}^{n}r_{s,u}\left\{ 
\begin{array}{c}
1,\text{ if }u=i; \\ 
-1,\text{ if }u=j; \\ 
0,\text{ if }u\neq i\text{ and }u\neq j%
\end{array}
\right. =r_{s,i}-r_{s,j},
\end{align*}
and thus%
\begin{align*}
g\left( e_{i}-e_{j}\right) & =\sum_{u=1}^{n}a_{u}\left\vert \left(
e_{i}-e_{j}\right) _{u}\right\vert +a\left\vert \left( e_{i}-e_{j}\right)
_{1}+\left( e_{i}-e_{j}\right) _{2}+...+\left( e_{i}-e_{j}\right)
_{n}\right\vert -\sum_{s\in S}b_{s}\left\vert r_{s}\left( e_{i}-e_{j}\right)
\right\vert \\
& =\sum_{u=1}^{n}a_{u}\left\vert \left\{ 
\begin{array}{c}
1,\text{ if }u=i; \\ 
-1,\text{ if }u=j; \\ 
0,\text{ if }u\neq i\text{ and }u\neq j%
\end{array}
\right. \right\vert +a\left\vert 0\right\vert -\sum_{s\in S}b_{s}\left\vert
r_{s,i}-r_{s,j}\right\vert \\
& =\left( a_{i}\left\vert 1\right\vert +a_{j}\left\vert -1\right\vert
\right) +a\left\vert 0\right\vert -\sum_{s\in S}b_{s}\left\vert
r_{s,i}-r_{s,j}\right\vert =\left( a_{i}+a_{j}\right) +0-\sum_{s\in
S}b_{s}\left\vert r_{s,i}-r_{s,j}\right\vert \\
& =a_{i}+a_{j}-\sum_{s\in S}b_{s}\left\vert r_{s,i}-r_{s,j}\right\vert \geq0
\end{align*}
(since $a_{i}+a_{j}\geq\sum\limits_{s\in S}b_{s}\left\vert
r_{s,i}-r_{s,j}\right\vert $ by the condition of Theorem 13).

So we have shown that $g\left( e_{i}\right) \geq0$ for every integer $%
i\in\left\{ 1,2,...,n\right\} ,$ and $g\left( e_{i}-e_{j}\right) \geq0$ for
any two distinct integers $i$ and $j$ from $\left\{ 1,2,...,n\right\} .$
Thus, Assertion $\mathcal{B}_{2}$ of Theorem 12 is fulfilled. According to
Theorem 12, the assertions $\mathcal{B}_{1}$ and $\mathcal{B}_{2}$ are
equivalent, so that Assertion $\mathcal{B}_{1}$ must be fulfilled as well.
Hence, $g\left( x\right) \geq0$ for every $x\in\mathbb{R}^{n}.$ In
particular, if we set $x=\left( 
\begin{array}{c}
y_{1} \\ 
y_{2} \\ 
... \\ 
y_{n}%
\end{array}
\right) ,$ then $r_{s}x=\sum\limits_{v=1}^{n}r_{s,v}y_{v}$ (since $%
r_{s}=\left( r_{s,1},r_{s,2},...,r_{s,n}\right) $), so that 
\begin{align*}
g\left( x\right) & =\sum_{u=1}^{n}a_{u}\left\vert y_{u}\right\vert
+a\left\vert y_{1}+y_{2}+...+y_{n}\right\vert -\sum_{s\in S}b_{s}\left\vert
r_{s}x\right\vert \\
& =\sum_{u=1}^{n}a_{u}\left\vert y_{u}\right\vert +a\left\vert
y_{1}+y_{2}+...+y_{n}\right\vert -\sum_{s\in S}b_{s}\left\vert
\sum\limits_{v=1}^{n}r_{s,v}y_{v}\right\vert \\
& =\sum\limits_{i=1}^{n}a_{i}\left\vert y_{i}\right\vert +a\left\vert
\sum\limits_{v=1}^{n}y_{v}\right\vert -\sum_{s\in S}b_{s}\left\vert
\sum\limits_{v=1}^{n}r_{s,v}y_{v}\right\vert ,
\end{align*}
and thus $g\left( x\right) \geq0$ yields%
\begin{equation*}
\sum\limits_{i=1}^{n}a_{i}\left\vert y_{i}\right\vert +a\left\vert
\sum\limits_{v=1}^{n}y_{v}\right\vert -\sum_{s\in S}b_{s}\left\vert
\sum\limits_{v=1}^{n}r_{s,v}y_{v}\right\vert \geq0.
\end{equation*}
Theorem 13 is thus proven.

\begin{center}
\textbf{6. A general condition for Popoviciu-like inequalities}
\end{center}

Now, we state a result more general than Theorem 5b:

\begin{quote}
\textbf{Theorem 14.} Let $n$ be a nonnegative integer. Let $a_{1},$ $a_{2},$ 
$...,$ $a_{n}$ and $a$ be $n+1$ nonnegative reals. Let $S$ be a finite set.
For every $s\in S,$ let $r_{s,1},$ $r_{s,2},$ $...,$ $r_{s,n}$ be $n$
nonnegative reals, and let $b_{s}$ be a nonnegative real. Assume that the
following two assertions hold\footnote{%
The second of these two assertions ($a_{i}+a_{j}\geq\sum\limits_{s\in
S}b_{s}\left\vert r_{s,i}-r_{s,j}\right\vert \ $for any two distinct
integers $i$ and $j$ from $\left\{ 1,2,...,n\right\} $) is identic with the
second required assertion in Theorem 13, but the first one ($%
a_{i}+a=\sum\limits_{s\in S}b_{s}r_{s,i}$ for every $i\in\left\{
1,2,...,n\right\} $) is \textit{stronger} than the first required assertion
in Theorem 13 (which only said that $a_{i}+a\geq \sum\limits_{s\in
S}b_{s}r_{s,i}$ for every $i\in\left\{ 1,2,...,n\right\} $).}:%
\begin{align*}
a_{i}+a & =\sum_{s\in S}b_{s}r_{s,i}\ \ \ \ \ \ \ \ \ \ \text{for every }%
i\in\left\{ 1,2,...,n\right\} ; \\
a_{i}+a_{j} & \geq\sum_{s\in S}b_{s}\left\vert r_{s,i}-r_{s,j}\right\vert \
\ \ \ \ \ \ \ \ \ \text{for any two distinct integers }i\text{ and }j\text{
from }\left\{ 1,2,...,n\right\} .
\end{align*}

Let $f$ be a convex function from an interval $I\subseteq\mathbb{R}$ to $%
\mathbb{R}.$ Let $w_{1},$ $w_{2},$ $...,$ $w_{n}$ be nonnegative reals.
Assume that $\sum\limits_{v=1}^{n}w_{v}\neq0$ and $\sum%
\limits_{v=1}^{n}r_{s,v}w_{v}\neq0$ for all $s\in S.$

Let $x_{1},$ $x_{2},$ $...,$ $x_{n}$ be $n$ points from the interval $I.$
Then, the inequality%
\begin{equation*}
\sum_{i=1}^{n}a_{i}w_{i}f\left( x_{i}\right) +a\left(
\sum_{v=1}^{n}w_{v}\right) f\left( \frac{\sum\limits_{v=1}^{n}w_{v}x_{v}}{%
\sum \limits_{v=1}^{n}w_{v}}\right) \geq\sum_{s\in S}b_{s}\left(
\sum_{v=1}^{n}r_{s,v}w_{v}\right) f\left( \frac{\sum%
\limits_{v=1}^{n}r_{s,v}w_{v}x_{v}}{\sum\limits_{v=1}^{n}r_{s,v}w_{v}}\right)
\end{equation*}
holds.

\textit{Remark.} Written in a less formal way, this inequality states that%
\begin{align*}
& \sum_{i=1}^{n}a_{i}w_{i}f\left( x_{i}\right) +a\left(
w_{1}+w_{2}+...+w_{n}\right) f\left( \frac{%
w_{1}x_{1}+w_{2}x_{2}+...+w_{n}x_{n}}{w_{1}+w_{2}+...+w_{n}}\right) \\
& \geq\sum_{s\in S}b_{s}\left(
r_{s,1}w_{1}+r_{s,2}w_{2}+...+r_{s,n}w_{n}\right) f\left( \frac{%
r_{s,1}w_{1}x_{1}+r_{s,2}w_{2}x_{2}+...+r_{s,n}w_{n}x_{n}}{%
r_{s,1}w_{1}+r_{s,2}w_{2}+...+r_{s,n}w_{n}}\right) .
\end{align*}
\end{quote}

\textit{Proof of Theorem 14.} Since the elements of the finite set $S$ are
used as labels only, we can assume without loss of generality that $%
S=\left\{ n+2,n+3,...,N\right\} $ for some integer $N\geq n+1$ (we just
rename the elements of $S$ into $n+2,$ $n+3,$ $...,$ $N,$ where $%
N=n+1+\left\vert S\right\vert ;$ this is possible because the set $S$ is
finite\footnote{%
In particular, $N=n+1$ if $S=\varnothing.$}). Define%
\begin{align*}
u_{i} & =a_{i}w_{i}\ \ \ \ \ \ \ \ \ \ \text{for all }i\in\left\{
1,2,...,n\right\} ; \\
u_{n+1} & =a\left( \sum_{v=1}^{n}w_{v}\right) ; \\
u_{s} & =-b_{s}\left( \sum_{v=1}^{n}r_{s,v}w_{v}\right) \ \ \ \ \ \ \ \ \ \ 
\text{for all }s\in\left\{ n+2,n+3,...,N\right\} \text{ (that is, for all }%
s\in S\text{).}
\end{align*}
Also define%
\begin{align*}
z_{i} & =x_{i}\ \ \ \ \ \ \ \ \ \ \text{for all }i\in\left\{
1,2,...,n\right\} ; \\
z_{n+1} & =\frac{\sum\limits_{v=1}^{n}w_{v}x_{v}}{\sum\limits_{v=1}^{n}w_{v}}%
; \\
z_{s} & =\frac{\sum\limits_{v=1}^{n}r_{s,v}w_{v}x_{v}}{\sum%
\limits_{v=1}^{n}r_{s,v}w_{v}}\ \ \ \ \ \ \ \ \ \ \text{for all }s\in\left\{
n+2,n+3,...,N\right\} \text{ (that is, for all }s\in S\text{).}
\end{align*}
Each of these $N$ reals $z_{1},$ $z_{2},$ $...,$ $z_{N}$ is a weighted mean
of the reals $x_{1},$ $x_{2},$ $...,$ $x_{n}$ with nonnegative weights.
Since the reals $x_{1},$ $x_{2},$ $...,$ $x_{n}$ lie in the interval $I,$ we
can thus conclude that each of the $N$ reals $z_{1},$ $z_{2},$ $...,$ $z_{N}$
lies in the interval $I$ as well. In other words, the points $z_{1},$ $%
z_{2}, $ $...,$ $z_{N}$ are $N$ points from $I.$

Now,%
\begin{align*}
& \sum_{i=1}^{n}a_{i}w_{i}f\left( x_{i}\right) +a\left(
\sum_{v=1}^{n}w_{v}\right) f\left( \frac{\sum\limits_{v=1}^{n}w_{v}x_{v}}{%
\sum\limits_{v=1}^{n}w_{v}}\right) -\sum_{s\in S}b_{s}\left(
\sum_{v=1}^{n}r_{s,v}w_{v}\right) f\left( \frac{\sum%
\limits_{v=1}^{n}r_{s,v}w_{v}x_{v}}{\sum\limits_{v=1}^{n}r_{s,v}w_{v}}\right)
\\
& =\sum_{i=1}^{n}\underbrace{a_{i}w_{i}}_{=u_{i}}f\left( \underbrace{x_{i}}%
_{=z_{i}}\right) +\underbrace{a\left( \sum_{v=1}^{n}w_{v}\right) }%
_{=u_{n+1}}f\left( \underbrace{\frac{\sum\limits_{v=1}^{n}w_{v}x_{v}}{%
\sum\limits_{v=1}^{n}w_{v}}}_{=z_{n+1}}\right) +\sum_{s\in S}\left( 
\underbrace{-b_{s}\left( \sum_{v=1}^{n}r_{s,v}w_{v}\right) }_{=u_{s}}\right)
f\left( \underbrace{\frac{\sum\limits_{v=1}^{n}r_{s,v}w_{v}x_{v}}{%
\sum\limits_{v=1}^{n}r_{s,v}w_{v}}}_{=z_{s}}\right) \\
& =\sum_{i=1}^{n}u_{i}f\left( z_{i}\right) +u_{n+1}f\left( z_{n+1}\right)
+\sum_{s\in S}u_{s}f\left( z_{s}\right) =\sum_{i=1}^{n}u_{i}f\left(
z_{i}\right) +u_{n+1}f\left( z_{n+1}\right) +\sum_{s=n+2}^{N}u_{s}f\left(
z_{s}\right) \\
& =\sum_{k=1}^{N}u_{k}f\left( z_{k}\right) .
\end{align*}
Hence, once we are able to show that $\sum\limits_{k=1}^{N}u_{k}f\left(
z_{k}\right) \geq0,$ we will obtain%
\begin{equation*}
\sum_{i=1}^{n}a_{i}w_{i}f\left( x_{i}\right) +a\left(
\sum_{v=1}^{n}w_{v}\right) f\left( \frac{\sum\limits_{v=1}^{n}w_{v}x_{v}}{%
\sum \limits_{v=1}^{n}w_{v}}\right) \geq\sum_{s\in S}b_{s}\left(
\sum_{v=1}^{n}r_{s,v}w_{v}\right) f\left( \frac{\sum%
\limits_{v=1}^{n}r_{s,v}w_{v}x_{v}}{\sum\limits_{v=1}^{n}r_{s,v}w_{v}}%
\right) ,
\end{equation*}
and thus Theorem 14 will be established.

Therefore, in order to prove Theorem 14, it remains to prove the inequality $%
\sum\limits_{k=1}^{N}u_{k}f\left( z_{k}\right) \geq0.$

We have%
\begin{align*}
\sum_{k=1}^{N}u_{k} &
=\sum_{i=1}^{n}u_{i}+u_{n+1}+\sum_{s=n+2}^{N}u_{s}=%
\sum_{i=1}^{n}u_{i}+u_{n+1}+\sum_{s\in S}u_{s} \\
& =\sum_{i=1}^{n}a_{i}w_{i}+a\left( \sum_{v=1}^{n}w_{v}\right) +\sum_{s\in
S}\left( -b_{s}\left( \sum_{v=1}^{n}r_{s,v}w_{v}\right) \right) \\
& =\sum_{i=1}^{n}a_{i}w_{i}+a\left( \sum_{i=1}^{n}w_{i}\right) +\sum_{s\in
S}\left( -b_{s}\left( \sum_{i=1}^{n}r_{s,i}w_{i}\right) \right) \\
& =\sum_{i=1}^{n}a_{i}w_{i}+\sum_{i=1}^{n}aw_{i}-\sum_{i=1}^{n}\sum_{s\in
S}b_{s}r_{s,i}w_{i}=\sum_{i=1}^{n}\left( a_{i}w_{i}+aw_{i}-\sum_{s\in
S}b_{s}r_{s,i}w_{i}\right) \\
& =\sum_{i=1}^{n}\left( a_{i}+a-\sum_{s\in S}b_{s}r_{s,i}\right) w_{i} \\
& =\sum_{i=1}^{n}0w_{i}\ \ \ \ \ \ \ \ \ \ \left( 
\begin{array}{c}
\text{since }a_{i}+a=\sum\limits_{s\in S}b_{s}r_{s,i}\text{ by an assumption
of Theorem 14,} \\ 
\text{and thus }a_{i}+a-\sum\limits_{s\in S}b_{s}r_{s,i}=0%
\end{array}
\right) \\
& =0.
\end{align*}

Next, we are going to prove that $\sum\limits_{k=1}^{N}u_{k}\left\vert
z_{k}-t\right\vert \geq0$ holds for every $t\in\left\{
z_{1},z_{2},...,z_{N}\right\} .$ In fact, let $t\in\left\{
z_{1},z_{2},...,z_{N}\right\} $ be arbitrary. Set $y_{i}=w_{i}\left(
x_{i}-t\right) $ for every $i\in\left\{ 1,2,...,n\right\} .$ Then, for all $%
i\in\left\{ 1,2,...,n\right\} ,$ we have $w_{i}\left( z_{i}-t\right)
=w_{i}\left( x_{i}-t\right) =y_{i}.$ Furthermore,%
\begin{equation*}
z_{n+1}-t=\frac{\sum\limits_{v=1}^{n}w_{v}x_{v}}{\sum\limits_{v=1}^{n}w_{v}}%
-t=\frac{\sum\limits_{v=1}^{n}w_{v}x_{v}-\sum\limits_{v=1}^{n}w_{v}\cdot t}{%
\sum\limits_{v=1}^{n}w_{v}}=\frac{\sum\limits_{v=1}^{n}w_{v}\left(
x_{v}-t\right) }{\sum\limits_{v=1}^{n}w_{v}}=\frac{\sum\limits_{v=1}^{n}y_{v}%
}{\sum\limits_{v=1}^{n}w_{v}}.
\end{equation*}
Finally, for all $s\in\left\{ n+2,n+3,...,N\right\} $ (that is, for all $%
s\in S$), we have%
\begin{equation*}
z_{s}-t=\frac{\sum\limits_{v=1}^{n}r_{s,v}w_{v}x_{v}}{\sum%
\limits_{v=1}^{n}r_{s,v}w_{v}}-t=\frac{\sum%
\limits_{v=1}^{n}r_{s,v}w_{v}x_{v}-\sum\limits_{v=1}^{n}r_{s,v}w_{v}\cdot t}{%
\sum\limits_{v=1}^{n}r_{s,v}w_{v}}=\frac{\sum\limits_{v=1}^{n}r_{s,v}w_{v}%
\left( x_{v}-t\right) }{\sum\limits_{v=1}^{n}r_{s,v}w_{v}}=\frac{%
\sum\limits_{v=1}^{n}r_{s,v}y_{v}}{\sum\limits_{v=1}^{n}r_{s,v}w_{v}}.
\end{equation*}
Hence,%
\begin{align*}
& \sum\limits_{k=1}^{N}u_{k}\left\vert z_{k}-t\right\vert =\sum
\limits_{i=1}^{n}u_{i}\left\vert z_{i}-t\right\vert +u_{n+1}\left\vert
z_{n+1}-t\right\vert +\sum\limits_{s=n+2}^{N}u_{s}\left\vert
z_{s}-t\right\vert \\
& =\sum\limits_{i=1}^{n}a_{i}\underbrace{w_{i}\left\vert z_{i}-t\right\vert }
_{\substack{ =\left\vert w_{i}\left( z_{i}-t\right) \right\vert ,  \\ \text{%
since }w_{i}\geq0}}+a\left( \sum_{v=1}^{n}w_{v}\right) \left\vert \frac{%
\sum\limits_{v=1}^{n}y_{v}}{\sum\limits_{v=1}^{n}w_{v}}\right\vert
+\sum\limits_{s=n+2}^{N}\left( -b_{s}\left(
\sum_{v=1}^{n}r_{s,v}w_{v}\right) \right) \left\vert \frac{%
\sum\limits_{v=1}^{n}r_{s,v}y_{v}}{\sum\limits_{v=1}^{n}r_{s,v}w_{v}}%
\right\vert \\
& =\sum\limits_{i=1}^{n}a_{i}\left\vert w_{i}\left( z_{i}-t\right)
\right\vert +a\left( \sum_{v=1}^{n}w_{v}\right) \frac{\left\vert
\sum\limits_{v=1}^{n}y_{v}\right\vert }{\sum\limits_{v=1}^{n}w_{v}}%
+\sum\limits_{s=n+2}^{N}\left( -b_{s}\left(
\sum_{v=1}^{n}r_{s,v}w_{v}\right) \right) \frac{\left\vert
\sum\limits_{v=1}^{n}r_{s,v}y_{v}\right\vert }{\sum%
\limits_{v=1}^{n}r_{s,v}w_{v}} \\
& \ \ \ \ \ \ \ \ \ \ \left( 
\begin{array}{c}
\text{here we have pulled the }\sum\limits_{v=1}^{n}w_{v}\text{ and }%
\sum\limits_{v=1}^{n}r_{s,v}w_{v}\text{ terms out of the modulus} \\ 
\text{signs, since they are positive (in fact, they are }\neq0\text{ by an
assumption} \\ 
\text{of Theorem 14, and nonnegative because }w_{i}\text{ and }r_{s,i}\text{
are all nonnegative)}%
\end{array}
\right) \\
& =\sum\limits_{i=1}^{n}a_{i}\left\vert y_{i}\right\vert +a\left\vert
\sum\limits_{v=1}^{n}y_{v}\right\vert +\sum\limits_{s=n+2}^{N}\left(
-b_{s}\right) \left\vert \sum\limits_{v=1}^{n}r_{s,v}y_{v}\right\vert
=\sum\limits_{i=1}^{n}a_{i}\left\vert y_{i}\right\vert +a\left\vert
\sum\limits_{v=1}^{n}y_{v}\right\vert
-\sum\limits_{s=n+2}^{N}b_{s}\left\vert
\sum\limits_{v=1}^{n}r_{s,v}y_{v}\right\vert \\
& =\sum\limits_{i=1}^{n}a_{i}\left\vert y_{i}\right\vert +a\left\vert
\sum\limits_{v=1}^{n}y_{v}\right\vert -\sum\limits_{s\in S}b_{s}\left\vert
\sum\limits_{v=1}^{n}r_{s,v}y_{v}\right\vert \geq0
\end{align*}
by Theorem 13 (in fact, we were allowed to apply Theorem 13 because all the
requirements of Theorem 13 are fulfilled - in particular, we have $%
a_{i}+a\geq\sum\limits_{s\in S}b_{s}r_{s,i}$ for every $i\in\left\{
1,2,...,n\right\} $ because we know that $a_{i}+a=\sum\limits_{s\in
S}b_{s}r_{s,i}$ for every $i\in\left\{ 1,2,...,n\right\} $ by an assumption
of Theorem 14).

Altogether, we have now shown the following: The points $z_{1},$ $z_{2},$ $%
...,$ $z_{N}$ are $N$ points from $I.$ The $N$ reals $u_{1},$ $u_{2},$ $...,$
$u_{N}$ satisfy $\sum\limits_{k=1}^{N}u_{k}=0,$ and $\sum%
\limits_{k=1}^{N}u_{k}\left\vert z_{k}-t\right\vert \geq 0$ holds for every $%
t\in \left\{ z_{1},z_{2},...,z_{N}\right\} .$ Hence, according to Theorem
8b, we have $\sum\limits_{k=1}^{N}u_{k}f\left( z_{k}\right) \geq 0.$ And as
we have seen above, once $\sum\limits_{k=1}^{N}u_{k}f\left( z_{k}\right)
\geq 0$ is shown, the proof of Theorem 14 is complete. Thus, Theorem 14 is
proven.

Theorem 14 gives a sufficient criterion for the validity of inequalities of
the kind%
\begin{align*}
& \text{convex combination of }f\left( x_{1}\right) ,\text{ }f\left(
x_{2}\right) ,\text{ }...,\text{ }f\left( x_{n}\right) \\
& \ \ \ \ \ \ \ \ \ \ \text{ and }f\left( \text{some weighted mean of }x_{1},%
\text{ }x_{2},\text{ }...,\text{ }x_{n}\right) \\
& \geq \text{convex combination of finitely many }f\left( \text{some other
weighted means of }x_{1},\text{ }x_{2},\text{ }...,\text{ }x_{n}\right) 
\text{'s,}
\end{align*}%
where $f$ is a convex function and $x_{1},$ $x_{2},$ $...,$ $x_{n}$ are $n$
reals in its domain, and where the weights of the weighted mean on the left
hand side are positive (those of the weighted means on the right hand side
may be $0$ as well, but still have to be nonnegative). This criterion turns
out to be necessary as well:

\begin{quote}
\textbf{Theorem 14b.} Let $n$ be a nonnegative integer. Let $w_{1},$ $w_{2},$
$...,$ $w_{n}$ be positive reals. Let $a_{1},$ $a_{2},$ $...,$ $a_{n}$ and $%
a $ be $n+1$ nonnegative reals. Let $S$ be a finite set. For every $s\in S,$
let $r_{s,1},$ $r_{s,2},$ $...,$ $r_{s,n}$ be $n$ nonnegative reals, and let 
$b_{s}$ be a nonnegative real. Let $I\subseteq \mathbb{R}$ be an interval.

Assume that the inequality%
\begin{equation*}
\sum_{i=1}^{n}a_{i}w_{i}f\left( x_{i}\right) +a\left(
\sum_{v=1}^{n}w_{v}\right) f\left( \frac{\sum\limits_{v=1}^{n}w_{v}x_{v}}{%
\sum\limits_{v=1}^{n}w_{v}}\right) \geq \sum_{s\in S}b_{s}\left(
\sum_{v=1}^{n}r_{s,v}w_{v}\right) f\left( \frac{\sum%
\limits_{v=1}^{n}r_{s,v}w_{v}x_{v}}{\sum\limits_{v=1}^{n}r_{s,v}w_{v}}\right)
\end{equation*}%
holds for any convex function $f:I\rightarrow \mathbb{R}$ and any $n$ points 
$x_{1},$ $x_{2},$ $...,$ $x_{n}$ in the interval $I.$ Then,%
\begin{align*}
a_{i}+a& =\sum_{s\in S}b_{s}r_{s,i}\ \ \ \ \ \ \ \ \ \ \text{for every }i\in
\left\{ 1,2,...,n\right\} ; \\
a_{i}+a_{j}& \geq \sum_{s\in S}b_{s}\left\vert r_{s,i}-r_{s,j}\right\vert \
\ \ \ \ \ \ \ \ \ \text{for any two distinct integers }i\text{ and }j\text{
from }\left\{ 1,2,...,n\right\} .
\end{align*}
\end{quote}

Since we are not going to use this fact, we are not proving it either, but
the idea of the proof is the following: Assume WLOG that $I=\left[ -1,1%
\right] .$ For every $i\in \left\{ 1,2,...,n\right\} ,$ you get $a_{i}+a\geq
\sum\limits_{s\in S}b_{s}r_{s,i}$ (by considering the convex function $%
f\left( x\right) =x$ and the points $x_{k}=\left\{ 
\begin{array}{c}
1,\text{ if }k=i; \\ 
0,\text{ if }k\neq i%
\end{array}%
\right. $) and $a_{i}+a\leq \sum\limits_{s\in S}b_{s}r_{s,i}$ (by
considering the convex function $f\left( x\right) =-x$ and the same points),
so that $a_{i}+a=\sum\limits_{s\in S}b_{s}r_{s,i}.$ For any two distinct
integers $i$ and $j$ from $\left\{ 1,2,...,n\right\} ,$ you get $%
a_{i}+a_{j}\geq \sum\limits_{s\in S}b_{s}\left\vert
r_{s,i}-r_{s,j}\right\vert $ (by considering the convex function $f\left(
x\right) =\left\vert x\right\vert $ and the points $x_{k}=\left\{ 
\begin{array}{c}
1,\text{ if }k=i; \\ 
-1,\text{ if }k=j; \\ 
0,\text{ if }k\neq i\text{ and }k\neq j%
\end{array}%
\right. $). This altogether proves Theorem 14b.

\begin{center}
\textbf{7. Proving the Popoviciu inequality}
\end{center}

Now we can finally step to the \textit{proof of Theorem 5b:}

We assume that $n\geq2,$ because all cases where $n<2$ (that is, $n=1$ or $%
n=0$) can be checked manually (and are uninteresting).

Let $a_{i}=\dbinom{n-2}{m-1}$ for every $i\in\left\{ 1,2,...,n\right\} .$
Let $a=\dbinom{n-2}{m-2}.$ These reals $a_{1},$ $a_{2},$ $...,$ $a_{n}$ and $%
a$ are all nonnegative (since $n\geq2$ yields $n-2\geq0$ and thus $\dbinom{%
n-2}{t}\geq0$ for all integers $t$).

Let $S=\left\{ s\subseteq\left\{ 1,2,...,n\right\} \mid\left\vert
s\right\vert =m\right\} ;$ that is, we denote by $S$ the set of all $m$%
-element subsets of the set $\left\{ 1,2,...,n\right\} .$ This set $S$ is
obviously finite.

For every $s\in S,$ define $n$ reals $r_{s,1},$ $r_{s,2},$ $...,$ $r_{s,n}$
as follows:%
\begin{equation*}
r_{s,i}=\left\{ 
\begin{array}{c}
1,\text{ if }i\in s; \\ 
0,\text{ if }i\notin s%
\end{array}
\right. \ \ \ \ \ \ \ \ \ \ \text{for every }i\in\left\{ 1,2,...,n\right\} .
\end{equation*}
Obviously, these reals $r_{s,1},$ $r_{s,2},$ $...,$ $r_{s,n}$ are all
nonnegative. Also, for every $s\in S,$ set $b_{s}=1;$ then, $b_{s}$ is a
nonnegative real as well.

For every $i\in\left\{ 1,2,...,n\right\} ,$ we have%
\begin{align*}
\sum_{s\in S}b_{s}r_{s,i} & =\sum_{s\in S}1r_{s,i}=\sum_{s\in
S}r_{s,i}=\sum_{s\in S}\left\{ 
\begin{array}{c}
1,\text{ if }i\in s; \\ 
0,\text{ if }i\notin s%
\end{array}
\right. =\sum_{\substack{ s\subseteq\left\{ 1,2,...,n\right\} ;  \\ %
\left\vert s\right\vert =m}}\left\{ 
\begin{array}{c}
1,\text{ if }i\in s; \\ 
0,\text{ if }i\notin s%
\end{array}
\right. \\
& =\left( \text{number of }m\text{-element subsets }s\text{ of the set }%
\left\{ 1,2,...,n\right\} \text{ that contain }i\right) \\
& =\dbinom{n-1}{m-1},
\end{align*}
so that%
\begin{align}
a_{i}+a & =\dbinom{n-2}{m-1}+\dbinom{n-2}{m-2}=\dbinom{n-1}{m-1}  \notag \\
& \ \ \ \ \ \ \ \ \ \ \left( \text{by the recurrence relation of the
binomial coefficients}\right)  \notag \\
& =\sum_{s\in S}b_{s}r_{s,i}.  \label{11}
\end{align}

For any two distinct integers $i$ and $j$ from $\left\{ 1,2,...,n\right\} ,$
we have%
\begin{align*}
& \sum_{s\in S}b_{s}\left\vert r_{s,i}-r_{s,j}\right\vert =\sum_{s\in
S}1\left\vert r_{s,i}-r_{s,j}\right\vert =\sum_{s\in S}\left\vert
r_{s,i}-r_{s,j}\right\vert \\
& =\sum_{s\in S}\left\vert \left\{ 
\begin{array}{c}
1,\text{ if }i\in s; \\ 
0,\text{ if }i\notin s%
\end{array}
\right. -\left\{ 
\begin{array}{c}
1,\text{ if }j\in s; \\ 
0,\text{ if }j\notin s%
\end{array}
\right. \right\vert =\sum_{s\in S}\left\{ 
\begin{array}{c}
0,\text{ if }i\in s\text{ and }j\in s; \\ 
1,\text{ if }i\in s\text{ and }j\notin s; \\ 
1,\text{ if }i\notin s\text{ and }j\in s; \\ 
0,\text{ if }i\notin s\text{ and }j\notin s%
\end{array}
\right. \\
& =\sum_{s\in S}\left( \left\{ 
\begin{array}{c}
1,\text{ if }i\in s\text{ and }j\notin s; \\ 
0\text{ otherwise}%
\end{array}
\right. +\left\{ 
\begin{array}{c}
1,\text{ if }i\notin s\text{ and }j\in s; \\ 
0\text{ otherwise}%
\end{array}
\right. \right) \\
& =\sum_{s\in S}\left\{ 
\begin{array}{c}
1,\text{ if }i\in s\text{ and }j\notin s; \\ 
0\text{ otherwise}%
\end{array}
\right. +\sum_{s\in S}\left\{ 
\begin{array}{c}
1,\text{ if }i\notin s\text{ and }j\in s; \\ 
0\text{ otherwise}%
\end{array}
\right. \\
& =\sum_{\substack{ s\subseteq\left\{ 1,2,...,n\right\} ;  \\ \left\vert
s\right\vert =m}}\left\{ 
\begin{array}{c}
1,\text{ if }i\in s\text{ and }j\notin s; \\ 
0\text{ otherwise}%
\end{array}
\right. +\sum_{\substack{ s\subseteq\left\{ 1,2,...,n\right\} ;  \\ %
\left\vert s\right\vert =m}}\left\{ 
\begin{array}{c}
1,\text{ if }i\notin s\text{ and }j\in s; \\ 
0\text{ otherwise}%
\end{array}
\right. \\
& =\left( \text{number of }m\text{-element subsets }s\text{ of the set }%
\left\{ 1,2,...,n\right\} \text{ that contain }i\text{ but not }j\right) \\
& +\left( \text{number of }m\text{-element subsets }s\text{ of the set }%
\left\{ 1,2,...,n\right\} \text{ that contain }j\text{ but not }i\right) \\
& =\dbinom{n-2}{m-1}+\dbinom{n-2}{m-1}=a_{i}+a_{j},
\end{align*}
so that%
\begin{equation}
a_{i}+a_{j}=\sum_{s\in S}b_{s}\left\vert r_{s,i}-r_{s,j}\right\vert .
\label{12}
\end{equation}
Also,%
\begin{equation}
\sum_{v=1}^{n}w_{v}=w_{1}+w_{2}+...+w_{n}\neq0  \label{13}
\end{equation}
(by an assumption of Theorem 5b).

The elements of $S$ are all the $m$-element subsets of $\left\{
1,2,...,n\right\} .$ Hence, to every element $s\in S$ uniquely correspond $m$
integers $i_{1},$ $i_{2},$ $...,$ $i_{m}$ satisfying $1\leq
i_{1}<i_{2}<...<i_{m}\leq n$ and $s=\left\{ i_{1},i_{2},...,i_{m}\right\} $
(these $m$ integers $i_{1},$ $i_{2},$ $...,$ $i_{m}$ are the $m$ elements of 
$s$ in increasing order). And conversely, any $m$ integers $i_{1},$ $i_{2},$ 
$...,$ $i_{m}$ satisfying $1\leq i_{1}<i_{2}<...<i_{m}\leq n$ can be
obtained this way - in fact, they correspond to the $m$-element set $%
s=\left\{ i_{1},i_{2},...,i_{m}\right\} \in S.$ Given an element $s\in S$
and the corresponding $m$ integers $i_{1},$ $i_{2},$ $...,$ $i_{m},$ we can
write%
\begin{align*}
\sum\limits_{v=1}^{n}r_{s,v}w_{v} & =\sum\limits_{v=1}^{n}\left\{ 
\begin{array}{c}
1,\text{ if }v\in s; \\ 
0,\text{ if }v\notin s%
\end{array}
\right. \cdot w_{v}=\sum\limits_{v\in s}w_{v}=\sum_{v\in\left\{
i_{1},i_{2},...,i_{m}\right\} }w_{v}=w_{i_{1}}+w_{i_{2}}+...+w_{i_{m}}; \\
\sum\limits_{v=1}^{n}r_{s,v}w_{v}x_{v} & =\sum\limits_{v=1}^{n}\left\{ 
\begin{array}{c}
1,\text{ if }v\in s; \\ 
0,\text{ if }v\notin s%
\end{array}
\right. \cdot w_{v}x_{v}=\sum\limits_{v\in s}w_{v}x_{v} \\
& =\sum_{v\in\left\{ i_{1},i_{2},...,i_{m}\right\}
}w_{v}x_{v}=w_{i_{1}}x_{i_{1}}+w_{i_{2}}x_{i_{2}}+...+w_{i_{m}}x_{i_{m}}.
\end{align*}
From this, we can conclude that 
\begin{equation}
\sum\limits_{v=1}^{n}r_{s,v}w_{v}\neq0\ \ \ \ \ \ \ \ \ \ \text{for every }%
s\in S  \label{14}
\end{equation}
(because $\sum%
\limits_{v=1}^{n}r_{s,v}w_{v}=w_{i_{1}}+w_{i_{2}}+...+w_{i_{m}},$ and $%
w_{i_{1}}+w_{i_{2}}+...+w_{i_{m}}\neq0$ by an assumption of Theorem 5b), and
we can also conclude that 
\begin{align}
& \sum_{s\in S}\left( \sum_{v=1}^{n}r_{s,v}w_{v}\right) f\left( \frac {%
\sum\limits_{v=1}^{n}r_{s,v}w_{v}x_{v}}{\sum\limits_{v=1}^{n}r_{s,v}w_{v}}%
\right)  \notag \\
& =\sum_{1\leq i_{1}<i_{2}<...<i_{m}\leq n}\left(
w_{i_{1}}+w_{i_{2}}+...+w_{i_{m}}\right) f\left( \frac{%
w_{i_{1}}x_{i_{1}}+w_{i_{2}}x_{i_{2}}+...+w_{i_{m}}x_{i_{m}}}{%
w_{i_{1}}+w_{i_{2}}+...+w_{i_{m}}}\right) .  \notag \\
&  \label{15}
\end{align}

Using the conditions of Theorem 5b and the relations (11), (12), (13) and
(14), we see that all conditions of Theorem 14 are fulfilled. Thus, we can
apply Theorem 14, and obtain%
\begin{equation*}
\sum_{i=1}^{n}a_{i}w_{i}f\left( x_{i}\right) +a\left(
\sum_{v=1}^{n}w_{v}\right) f\left( \frac{\sum\limits_{v=1}^{n}w_{v}x_{v}}{%
\sum \limits_{v=1}^{n}w_{v}}\right) \geq\sum_{s\in S}b_{s}\left(
\sum_{v=1}^{n}r_{s,v}w_{v}\right) f\left( \frac{\sum%
\limits_{v=1}^{n}r_{s,v}w_{v}x_{v}}{\sum\limits_{v=1}^{n}r_{s,v}w_{v}}%
\right) .
\end{equation*}
This rewrites as%
\begin{align*}
& \sum_{i=1}^{n}\dbinom{n-2}{m-1}w_{i}f\left( x_{i}\right) +\dbinom {n-2}{m-2%
}\left( \sum_{v=1}^{n}w_{v}\right) f\left( \frac{\sum
\limits_{v=1}^{n}w_{v}x_{v}}{\sum\limits_{v=1}^{n}w_{v}}\right) \\
& \geq\sum_{s\in S}1\left( \sum_{v=1}^{n}r_{s,v}w_{v}\right) f\left( \frac{%
\sum\limits_{v=1}^{n}r_{s,v}w_{v}x_{v}}{\sum\limits_{v=1}^{n}r_{s,v}w_{v}}%
\right) .
\end{align*}
In other words,%
\begin{align*}
& \dbinom{n-2}{m-1}\sum_{i=1}^{n}w_{i}f\left( x_{i}\right) +\dbinom {n-2}{m-2%
}\left( \sum_{v=1}^{n}w_{v}\right) f\left( \frac{\sum
\limits_{v=1}^{n}w_{v}x_{v}}{\sum\limits_{v=1}^{n}w_{v}}\right) \\
& \geq\sum_{s\in S}\left( \sum_{v=1}^{n}r_{s,v}w_{v}\right) f\left( \frac{%
\sum\limits_{v=1}^{n}r_{s,v}w_{v}x_{v}}{\sum\limits_{v=1}^{n}r_{s,v}w_{v}}%
\right) .
\end{align*}
Using (15) and the obvious relations%
\begin{align*}
\sum_{v=1}^{n}w_{v} & =w_{1}+w_{2}+...+w_{n}; \\
\sum\limits_{v=1}^{n}w_{v}x_{v} & =w_{1}x_{1}+w_{2}x_{2}+...+w_{n}x_{n},
\end{align*}
we can rewrite this as%
\begin{align*}
& \dbinom{n-2}{m-1}\sum_{i=1}^{n}w_{i}f\left( x_{i}\right) +\dbinom {n-2}{m-2%
}\left( w_{1}+w_{2}+...+w_{n}\right) f\left( \frac{%
w_{1}x_{1}+w_{2}x_{2}+...+w_{n}x_{n}}{w_{1}+w_{2}+...+w_{n}}\right) \\
& \geq\sum_{1\leq i_{1}<i_{2}<...<i_{m}\leq n}\left(
w_{i_{1}}+w_{i_{2}}+...+w_{i_{m}}\right) f\left( \frac{%
w_{i_{1}}x_{i_{1}}+w_{i_{2}}x_{i_{2}}+...+w_{i_{m}}x_{i_{m}}}{%
w_{i_{1}}+w_{i_{2}}+...+w_{i_{m}}}\right) .
\end{align*}
This proves Theorem 5b.

\begin{center}
\textbf{8. A cyclic inequality}
\end{center}

The most general form of the Popoviciu inequality is now proven. But this is
not the end to the applications of Theorem 14. We will now apply it to show
a cyclic inequality similar to Popoviciu's:

\begin{quote}
\textbf{Theorem 15a.} Let $f$ be a convex function from an interval $%
I\subseteq\mathbb{R}$ to $\mathbb{R}.$ Let $x_{1},$ $x_{2},$ $...,$ $x_{n}$
be finitely many points from $I.$

We extend the indices in $x_{1},$ $x_{2},$ $...,$ $x_{n}$ cyclically modulo $%
n$; this means that for any integer $i\notin\left\{ 1,2,...,n\right\} ,$ we
define a real $x_{i}$ by setting $x_{i}=x_{j},$ where $j$ is the integer
from the set $\left\{ 1,2,...,n\right\} $ such that $i\equiv j\func{mod}n.$
(For instance, this means that $x_{n+3}=x_{3}.$)

Let $x=\dfrac{x_{1}+x_{2}+...+x_{n}}{n}$. Let $r$ be an integer. Then,%
\begin{equation*}
2\sum_{i=1}^{n}f\left( x_{i}\right) +n\left( n-2\right) f\left( x\right)
\geq n\sum_{s=1}^{n}f\left( x+\dfrac{x_{s}-x_{s+r}}{n}\right) .
\end{equation*}
\end{quote}

A weighted version of this inequality is:

\begin{quote}
\textbf{Theorem 15b.} Let $f$ be a convex function from an interval $%
I\subseteq\mathbb{R}$ to $\mathbb{R}.$ Let $x_{1},$ $x_{2},$ $...,$ $x_{n}$
be finitely many points from $I.$ Let $r$ be an integer.

Let $w_{1},$ $w_{2},$ $...,$ $w_{n}$ be nonnegative reals. Let $x=\dfrac {%
\sum\limits_{v=1}^{n}w_{v}x_{v}}{\sum\limits_{v=1}^{n}w_{v}}$ and $%
w=\sum\limits_{v=1}^{n}w_{v}$. Assume that $w\neq0$ and that $w+\left(
w_{s}-w_{s+r}\right) \neq0$ for every $s\in S.$

We extend the indices in $x_{1},$ $x_{2},$ $...,$ $x_{n}$ and in $w_{1},$ $%
w_{2},$ $...,$ $w_{n}$ cyclically modulo $n$; this means that for any
integer $i\notin\left\{ 1,2,...,n\right\} ,$ we define reals $x_{i}$ and $%
w_{i}$ by setting $x_{i}=x_{j}$ and $w_{i}=w_{j},$ where $j$ is the integer
from the set $\left\{ 1,2,...,n\right\} $ such that $i\equiv j\func{mod}n.$
(For instance, this means that $x_{n+3}=x_{3}$ and $w_{n+2}=w_{2}.$)

Then,%
\begin{equation*}
2\sum_{i=1}^{n}w_{i}f\left( x_{i}\right) +\left( n-2\right) wf\left(
x\right) \geq\sum_{s=1}^{n}\left( w+\left( w_{s}-w_{s+r}\right) \right)
f\left( \frac{\sum\limits_{v=1}^{n}w_{v}x_{v}+\left(
w_{s}x_{s}-w_{s+r}x_{s+r}\right) }{w+\left( w_{s}-w_{s+r}\right) }\right) .
\end{equation*}
\end{quote}

\textit{Proof of Theorem 15b.} We assume that $n\geq2,$ because all cases
where $n<2$ (that is, $n=1$ or $n=0$) can be checked manually (and are
uninteresting).

Before we continue with the proof, let us introduce a simple notation: For
any assertion $\mathcal{A}$, we denote by $\left[ \mathcal{A}\right] $ the
Boolean value of the assertion $\mathcal{A}$ (that is, $\left[ \mathcal{A}%
\right] =\left\{ 
\begin{array}{c}
1\text{, if }\mathcal{A}\text{ is true;} \\ 
0\text{, if }\mathcal{A}\text{ is false}%
\end{array}
\right. $). Therefore, $0\leq\left[ \mathcal{A}\right] \leq1$ for every
assertion $\mathcal{A}.$

Let $a_{i}=2$ for every $i\in\left\{ 1,2,...,n\right\} .$ Let $a=n-2.$ These
reals $a_{1},$ $a_{2},$ $...,$ $a_{n}$ and $a$ are all nonnegative (since $%
n\geq2$ yields $n-2\geq0$).

Let $S=\left\{ 1,2,...,n\right\} .$ This set $S$ is obviously finite.

For every $s\in S,$ define $n$ reals $r_{s,1},$ $r_{s,2},$ $...,$ $r_{s,n}$
as follows:%
\begin{equation*}
r_{s,i}=1+\left[ i=s\right] -\left[ i\equiv s+r\func{mod}n\right] \ \ \ \ \
\ \ \ \ \ \text{for every }i\in\left\{ 1,2,...,n\right\} .
\end{equation*}
These reals $r_{s,1},$ $r_{s,2},$ $...,$ $r_{s,n}$ are all nonnegative
(because%
\begin{equation*}
r_{s,i}=1+\underbrace{\left[ i=s\right] }_{\geq0}-\underbrace{\left[ i\equiv
s+r\func{mod}n\right] }_{\leq1}\geq1+0-1=0
\end{equation*}
for every $i\in\left\{ 1,2,...,n\right\} $). Also, for every $s\in S,$ set $%
b_{s}=1;$ then, $b_{s}$ is a nonnegative real as well.

For every $i\in\left\{ 1,2,...,n\right\} ,$ we have $\sum\limits_{s=1}^{n}%
\left[ i=s\right] =1$ (because there exists one and only one $s\in\left\{
1,2,...,n\right\} $ satisfying $i=s$). Also, for every $i\in\left\{
1,2,...,n\right\} ,$ we have $\sum\limits_{s=1}^{n}\left[ s\equiv i-r\func{%
mod}n\right] =1$ (because there exists one and only one $s\in\left\{
1,2,...,n\right\} $ satisfying $s\equiv i-r\func{mod}n$). In other words, $%
\sum\limits_{s=1}^{n}\left[ i\equiv s+r\func{mod}n\right] =1$ (because $%
\left[ s\equiv i-r\func{mod}n\right] =\left[ i\equiv s+r\func{mod}n\right] $%
).

For every $i\in\left\{ 1,2,...,n\right\} ,$ we have%
\begin{align*}
\sum_{s\in S}b_{s}r_{s,i} & =\sum_{s=1}^{n}\underbrace{b_{s}}%
_{=1}r_{s,i}=\sum_{s=1}^{n}r_{s,i}=\sum_{s=1}^{n}\left( 1+\left[ i=s\right] -%
\left[ i\equiv s+r\func{mod}n\right] \right) \\
& =\sum_{s=1}^{n}1+\sum_{s=1}^{n}\left[ i=s\right] -\sum_{s=1}^{n}\left[
i\equiv s+r\func{mod}n\right] =n+1-1=n=2+\left( n-2\right) =a_{i}+a,
\end{align*}
so that%
\begin{equation}
a_{i}+a=\sum_{s\in S}b_{s}r_{s,i}.  \label{16}
\end{equation}

For any two integers $i$ and $j$ from $\left\{ 1,2,...,n\right\} ,$ we have%
\begin{align*}
\sum_{s=1}^{n}\left\vert r_{s,i}-1\right\vert & =\sum_{s=1}^{n}\left\vert
\left( 1+\left[ i=s\right] -\left[ i\equiv s+r\func{mod}n\right] \right)
-1\right\vert \\
& =\sum_{s=1}^{n}\left\vert \left[ i=s\right] +\left( -\left[ i\equiv s+r%
\func{mod}n\right] \right) \right\vert \leq\sum_{s=1}^{n}\left( \left\vert %
\left[ i=s\right] \right\vert +\left\vert -\left[ i\equiv s+r\func{mod}n%
\right] \right\vert \right) \\
& \ \ \ \ \ \ \ \ \ \ \left( \text{by the triangle inequality}\right) \\
& =\sum_{s=1}^{n}\left( \left[ i=s\right] +\left[ i\equiv s+r\func{mod}n%
\right] \right) \\
& \ \ \ \ \ \ \ \ \ \ \left( 
\begin{array}{c}
\text{because }\left[ i=s\right] \text{ and }\left[ i\equiv s+r\func{mod}n%
\right] \text{ are nonnegative, so that} \\ 
\left\vert \left[ i=s\right] \right\vert =\left[ i=s\right] \text{ and }%
\left\vert -\left[ i\equiv s+r\func{mod}n\right] \right\vert =\left[ i\equiv
s+r\func{mod}n\right]%
\end{array}
\right) \\
& =\sum_{s=1}^{n}\left[ i=s\right] +\sum_{s=1}^{n}\left[ i\equiv s+r\func{mod%
}n\right] =1+1=2
\end{align*}
and similarly $\sum\limits_{s=1}^{n}\left\vert r_{s,j}-1\right\vert \leq2,$
so that%
\begin{align*}
\sum_{s\in S}b_{s}\left\vert r_{s,i}-r_{s,j}\right\vert & =\sum_{s=1}^{n}%
\underbrace{b_{s}}_{=1}\left\vert r_{s,i}-r_{s,j}\right\vert =\sum
_{s=1}^{n}\left\vert r_{s,i}-r_{s,j}\right\vert =\sum_{s=1}^{n}\left\vert
\left( r_{s,i}-1\right) +\left( 1-r_{s,j}\right) \right\vert \\
& \leq\sum_{s=1}^{n}\left( \left\vert r_{s,i}-1\right\vert +\left\vert
1-r_{s,j}\right\vert \right) \ \ \ \ \ \ \ \ \ \ \left( \text{by the
triangle inequality}\right) \\
& =\sum_{s=1}^{n}\left( \left\vert r_{s,i}-1\right\vert +\left\vert
r_{s,j}-1\right\vert \right) =\sum_{s=1}^{n}\left\vert r_{s,i}-1\right\vert
+\sum_{s=1}^{n}\left\vert r_{s,j}-1\right\vert \\
& \leq2+2=a_{i}+a_{j},
\end{align*}
and thus%
\begin{equation}
a_{i}+a_{j}\geq\sum_{s\in S}b_{s}\left\vert r_{s,i}-r_{s,j}\right\vert .
\label{17}
\end{equation}

For every $s\in S$ (that is, for every $s\in\left\{ 1,2,...,n\right\} $), we
have%
\begin{equation*}
\sum\limits_{v=1}^{n}\left[ v\equiv s+r\func{mod}n\right] \cdot
w_{v}=\sum\limits_{v=1}^{n}\left\{ 
\begin{array}{c}
w_{v},\text{ if }v\equiv s+r\func{mod}n; \\ 
0\text{ otherwise}%
\end{array}
\right. =w_{s+r}
\end{equation*}
(because there is one and only one element $v\in\left\{ 1,2,...,n\right\} $
that satisfies $v\equiv s+r\func{mod}n,$ and for this element $v,$ we have $%
w_{v}=w_{s+r}$), so that%
\begin{align*}
\sum\limits_{v=1}^{n}r_{s,v}w_{v} & =\sum\limits_{v=1}^{n}\left( 1+\left[ v=s%
\right] -\left[ v\equiv s+r\func{mod}n\right] \right) \cdot w_{v} \\
& =\underbrace{\sum\limits_{v=1}^{n}w_{v}}_{=w}+\underbrace{\sum
\limits_{v=1}^{n}\left[ v=s\right] \cdot w_{v}}_{=w_{s}}-\underbrace {%
\sum\limits_{v=1}^{n}\left[ v\equiv s+r\func{mod}n\right] \cdot w_{v}}%
_{=w_{s+r}} \\
& =w+w_{s}-w_{s+r}=w+\left( w_{s}-w_{s+r}\right) .
\end{align*}
Also, for every $s\in S$ (that is, for every $s\in\left\{ 1,2,...,n\right\} $%
), we have%
\begin{equation*}
\sum\limits_{v=1}^{n}\left[ v\equiv s+r\func{mod}n\right] \cdot
w_{v}x_{v}=\sum\limits_{v=1}^{n}\left\{ 
\begin{array}{c}
w_{v}x_{v},\text{ if }v\equiv s+r\func{mod}n; \\ 
0\text{ otherwise}%
\end{array}
\right. =w_{s+r}x_{s+r}
\end{equation*}
(because there is one and only one element $v\in\left\{ 1,2,...,n\right\} $
that satisfies $v\equiv s+r\func{mod}n,$ and for this element $v,$ we have $%
w_{v}=w_{s+r}$ and $x_{v}=x_{s+r}$), and thus%
\begin{align*}
\sum\limits_{v=1}^{n}r_{s,v}w_{v}x_{v} & =\sum\limits_{v=1}^{n}\left( 1+ 
\left[ v=s\right] -\left[ v\equiv s+r\func{mod}n\right] \right) \cdot
w_{v}x_{v} \\
& =\sum\limits_{v=1}^{n}w_{v}x_{v}+\underbrace{\sum\limits_{v=1}^{n}\left[
v=s\right] \cdot w_{v}x_{v}}_{=w_{s}x_{s}}-\underbrace{\sum\limits_{v=1}^{n}%
\left[ v\equiv s+r\func{mod}n\right] \cdot w_{v}x_{v}}_{=w_{s+r}x_{s+r}} \\
& =\sum\limits_{v=1}^{n}w_{v}x_{v}+w_{s}x_{s}-w_{s+r}x_{s+r}=\sum
\limits_{v=1}^{n}w_{v}x_{v}+\left( w_{s}x_{s}-w_{s+r}x_{s+r}\right) .
\end{align*}

Now it is clear that $\sum\limits_{v=1}^{n}r_{s,v}w_{v}\neq0$ for all $s\in
S $ (because $\sum\limits_{v=1}^{n}r_{s,v}w_{v}=w+\left(
w_{s}-w_{s+r}\right) $ and $w+\left( w_{s}-w_{s+r}\right) \neq0$). Also, $%
\sum\limits_{v=1}^{n}w_{v}\neq0$ (since $\sum\limits_{v=1}^{n}w_{v}=w$ and $%
w\neq0$). Using these two relations, the conditions of Theorem 15b and the
relations (16) and (17), we see that all conditions of Theorem 14 are
fulfilled. Hence, we can apply Theorem 14 and obtain%
\begin{equation*}
\sum_{i=1}^{n}\underbrace{a_{i}}_{=2}w_{i}f\left( x_{i}\right) +\underbrace{a%
}_{=n-2}\left( \underbrace{\sum_{v=1}^{n}w_{v}}_{=w}\right) f\left( 
\underbrace{\frac{\sum\limits_{v=1}^{n}w_{v}x_{v}}{\sum
\limits_{v=1}^{n}w_{v}}}_{=x}\right) \geq\sum_{s\in S}\underbrace{b_{s}}%
_{=1}\left( \sum_{v=1}^{n}r_{s,v}w_{v}\right) f\left( \frac{\sum
\limits_{v=1}^{n}r_{s,v}w_{v}x_{v}}{\sum\limits_{v=1}^{n}r_{s,v}w_{v}}%
\right) .
\end{equation*}
This immediately simplifies to%
\begin{equation*}
\sum_{i=1}^{n}2w_{i}f\left( x_{i}\right) +\left( n-2\right) wf\left(
x\right) \geq\sum_{s\in S}1\left( \sum_{v=1}^{n}r_{s,v}w_{v}\right) f\left( 
\frac{\sum\limits_{v=1}^{n}r_{s,v}w_{v}x_{v}}{\sum%
\limits_{v=1}^{n}r_{s,v}w_{v}}\right) .
\end{equation*}
Recalling that for every $s\in S,$ we have $\sum%
\limits_{v=1}^{n}r_{s,v}w_{v}=w+\left( w_{s}-w_{s+r}\right) $ and $%
\sum\limits_{v=1}^{n}r_{s,v}w_{v}x_{v}=\sum\limits_{v=1}^{n}w_{v}x_{v}+%
\left( w_{s}x_{s}-w_{s+r}x_{s+r}\right) ,$ we can rewrite this as%
\begin{equation*}
\sum_{i=1}^{n}2w_{i}f\left( x_{i}\right) +\left( n-2\right) wf\left(
x\right) \geq\sum_{s\in S}1\left( w+\left( w_{s}-w_{s+r}\right) \right)
f\left( \frac{\sum\limits_{v=1}^{n}w_{v}x_{v}+\left(
w_{s}x_{s}-w_{s+r}x_{s+r}\right) }{w+\left( w_{s}-w_{s+r}\right) }\right) .
\end{equation*}
In other words,%
\begin{equation*}
2\sum_{i=1}^{n}w_{i}f\left( x_{i}\right) +\left( n-2\right) wf\left(
x\right) \geq\sum_{s=1}^{n}\left( w+\left( w_{s}-w_{s+r}\right) \right)
f\left( \frac{\sum\limits_{v=1}^{n}w_{v}x_{v}+\left(
w_{s}x_{s}-w_{s+r}x_{s+r}\right) }{w+\left( w_{s}-w_{s+r}\right) }\right) .
\end{equation*}
This proves Theorem 15b.

Theorem 15a follows from Theorem 15b if we set $w_{1}=w_{2}=...=w_{n}=1.$

Theorem 15a generalizes two inequalities that have appeared on the MathLinks
forum. The first of these results if we apply Theorem 15a to $r=1,$ to $r=2,$
to $r=3,$ and so on up to $r=n-1,$ and sum up the $n-1$ inequalities
obtained:

\begin{quote}
\textbf{Theorem 16.} Let $f$ be a convex function from an interval $%
I\subseteq\mathbb{R}$ to $\mathbb{R}.$ Let $x_{1},$ $x_{2},$ $...,$ $x_{n}$
be finitely many points from $I.$

Let $x=\dfrac{x_{1}+x_{2}+...+x_{n}}{n}$. Then,%
\begin{equation*}
2\left( n-1\right) \sum_{i=1}^{n}f\left( x_{i}\right) +n\left( n-1\right)
\left( n-2\right) f\left( x\right) \geq n\sum_{\substack{ 1\leq i\leq n;  \\ %
1\leq j\leq n;  \\ i\neq j}}f\left( x+\dfrac{x_{i}-x_{j}}{n}\right) .
\end{equation*}
\end{quote}

This inequality occured in [6], post \#4 as a result by Vasile C\^{\i}rtoaje
(Vasc). Our Theorem 15a is therefore a strengthening of this result.

The next inequality was proposed by Michael Rozenberg (aka "Arqady") in [7]:

\begin{quote}
\textbf{Theorem 17.} Let $a,$ $b,$ $c,$ $d$ be four nonnegative reals. Then,%
\begin{equation*}
a^{4}+b^{4}+c^{4}+d^{4}+4abcd\geq2\left(
a^{2}bc+b^{2}cd+c^{2}da+d^{2}ab\right) .
\end{equation*}
\end{quote}

\textit{Proof of Theorem 17.} The case when at least one of the reals $a,$ $%
b,$ $c,$ $d$ equals $0$ is easy (and a limiting case). Hence, we can assume
for the rest of this proof that none of the reals $a,$ $b,$ $c,$ $d$ equals $%
0.$ Since $a,$ $b,$ $c,$ $d$ are nonnegative, this means that $a,$ $b,$ $c,$ 
$d$ are positive.

Let $A=\ln\left( a^{4}\right) ,$ $B=\ln\left( b^{4}\right) ,$ $C=\ln\left(
c^{4}\right) ,$ $D=\ln\left( d^{4}\right) .$ Then, $\exp A=a^{4},$ $\exp
B=b^{4},$ $\exp C=c^{4},$ $\exp D=d^{4}.$

Let $I\subseteq\mathbb{R}$ be an interval containing the reals $A,$ $B,$ $C, 
$ $D$ (for instance, $I=\mathbb{R}$). Let $f:I\rightarrow\mathbb{R}$ be the
function defined by $f\left( x\right) =\exp x$ for all $x\in I.$ Then, it is
known that this function $f$ is convex. Applying Theorem 15a to $n=4,$ $%
x_{1}=A,$ $x_{2}=B,$ $x_{3}=C,$ $x_{4}=D,$ and $r=3,$ we obtain%
\begin{align*}
& 2\left( f\left( A\right) +f\left( B\right) +f\left( C\right) +f\left(
D\right) \right) +4\left( 4-2\right) f\left( \dfrac{A+B+C+D}{4}\right) \\
& \geq4\left( f\left( \dfrac{A+B+C+D}{4}+\dfrac{A-D}{4}\right) +f\left( 
\dfrac{A+B+C+D}{4}+\dfrac{B-A}{4}\right) \right. \\
& \ \ \ \ \ \ \ \ \ \ \left. +f\left( \dfrac{A+B+C+D}{4}+\dfrac{C-B}{4}%
\right) +f\left( \dfrac{A+B+C+D}{4}+\dfrac{D-C}{4}\right) \right) .
\end{align*}
Dividing this by $2$ and simplifying, we obtain%
\begin{align*}
& f\left( A\right) +f\left( B\right) +f\left( C\right) +f\left( D\right)
+4f\left( \dfrac{A+B+C+D}{4}\right) \\
& \geq2\left( f\left( \dfrac{2A+B+C}{4}\right) +f\left( \dfrac{2B+C+D}{4}%
\right) +f\left( \dfrac{2C+D+A}{4}\right) +f\left( \dfrac{2D+A+B}{4}\right)
\right) .
\end{align*}
Since we have%
\begin{align*}
f\left( A\right) & =\exp A=a^{4}\ \ \ \ \ \ \ \ \ \ \text{and similarly} \\
f\left( B\right) & =b^{4},\text{ }f\left( C\right) =c^{4},\text{ and }%
f\left( D\right) =d^{4}; \\
f\left( \dfrac{A+B+C+D}{4}\right) & =\exp\dfrac{A+B+C+D}{4}=\sqrt[4]{\exp
A\cdot\exp B\cdot\exp C\cdot\exp D} \\
& =\sqrt[4]{a^{4}\cdot b^{4}\cdot c^{4}\cdot d^{4}}=abcd; \\
f\left( \dfrac{2A+B+C}{4}\right) & =\exp\dfrac{2A+B+C}{4}=\sqrt[4]{\left(
\exp A\right) ^{2}\cdot\exp B\cdot\exp C} \\
& =\sqrt[4]{\left( a^{4}\right) ^{2}\cdot b^{4}\cdot c^{4}}=a^{2}bc\ \ \ \ \
\ \ \ \ \ \text{and similarly} \\
f\left( \dfrac{2B+C+D}{4}\right) & =b^{2}cd,\text{ }f\left( \dfrac{2C+D+A}{4}%
\right) =c^{2}da,\text{ and }f\left( \dfrac{2D+A+B}{4}\right) =d^{2}ab,
\end{align*}
this becomes%
\begin{equation*}
a^{4}+b^{4}+c^{4}+d^{4}+4abcd\geq2\left(
a^{2}bc+b^{2}cd+c^{2}da+d^{2}ab\right) .
\end{equation*}
This proves Theorem 17.

\begin{center}
\textbf{References}
\end{center}

[1] Vasile C\^{\i}rtoaje, \textit{Two Generalizations of Popoviciu's
Inequality}, Crux Mathematicorum 5/2001 (volume 31), pp. 313-318.\newline
\texttt{http://journals.cms.math.ca/CRUX/}

[2] Billzhao et al., \textit{Generalized Popoviciu - MathLinks topic \#19097}%
.\newline
\texttt{http://www.mathlinks.ro/Forum/viewtopic.php?t=19097}

[3] Billzhao et al., \textit{Like Popoviciu - MathLinks topic \#21786}.%
\newline
\texttt{http://www.mathlinks.ro/Forum/viewtopic.php?t=21786}

[4] Darij Grinberg et al., \textit{The Karamata Inequality - MathLinks topic
\#14975}.\newline
\texttt{http://www.mathlinks.ro/Forum/viewtopic.php?t=14975}

[5] Darij Grinberg et al., \textit{Weighted majorization and a result
stronger than Fuchs - MathLinks topic \#104714}.\newline
\texttt{http://www.mathlinks.ro/Forum/viewtopic.php?t=104714}

[6] Harazi et al., \textit{improvement of Popoviciu's inequality in a
particular case - MathLinks topic \#22364}.\newline
\texttt{http://www.mathlinks.ro/Forum/viewtopic.php?t=22364}

[7] Arqady et al., \textit{New, old inequality - MathLinks topic \#56040}.%
\newline
\texttt{http://www.mathlinks.ro/Forum/viewtopic.php?t=56040}

\end{document}